\documentclass[preprint,12pt]{elsarticle}

\makeatletter
\def\ps@pprintTitle{%
	\let\@oddhead\@empty
	\let\@evenhead\@empty
	\def\@oddfoot{\centerline{\thepage}}%
	\let\@evenfoot\@oddfoot}
\makeatother

\usepackage{amssymb}

\usepackage{xcolor}

\usepackage[colorlinks=true,linkcolor=blue,filecolor=blue,urlcolor=blue,citecolor=blue,pdftex,plainpages=false]{hyperref}

\begin{document}

\begin{frontmatter}



\title{Commutator-free Lie group methods with minimum storage requirements and reuse of exponentials}


\author{Alexei Bazavov}
\ead{bazavov@msu.edu}
\address{Department of Computational Mathematics, Science and Engineering and\\
Department of Physics and Astronomy,\\
Michigan State University, East Lansing, MI 48824, USA}

\begin{abstract}
A new format for commutator-free Lie group methods is proposed
based on explicit classical Runge-Kutta schemes.
In this format exponentials are reused at every stage and the storage is 
required only for two quantities: the right hand side of the differential equation
evaluated at a given Runge-Kutta stage and the function value updated at the
same stage. The next stage of the scheme is able to overwrite these values.
The result is proven for a 3-stage third order method and a conjecture for
higher order methods is formulated. Five numerical examples are provided in
support of the conjecture.
This new class of structure-preserving
integrators has a wide variety of applications for numerically
solving differential equations on manifolds.
\end{abstract}


\end{frontmatter}


\section{Introduction}

In many scientific and engineering applications there is a need to solve 
ordinary or partial differential equations numerically. A variety of methods exist
and one of the popular ones is the Runge-Kutta method~\cite{ButcherBook,HairerBook1}.
Often, one would like to build numerical schemes that preserve the \textit{structure}
of the original differential equations. For instance, for free rigid body rotation
the (properly normalized) vector of the angular momentum evolves on the $S^2$
manifold, \textit{i.e.} the surface of a three-dimensional sphere. It is beneficial
when a time-stepping scheme maintains this property at every step of the integration.

Ideas along these lines have been pursued over last three decades and lead to
development of geometric integrators~\cite{Hairer2006}, see
also~\cite{christiansen_munthe-kaas_owren_2011,CELLEDONI20141040} for
recent reviews. As argued in \cite{christiansen_munthe-kaas_owren_2011},
preservation of geometric properties is beneficial and often leads to increased
stability, smaller local error as well as slower global error growth in
long-time simulations. Many applications involve differential equation on Lie groups
or manifolds with Lie group action. The first major step in building Lie group
methods based on classical Runge-Kutta schemes was taken by 
Crouch and Grossman~\cite{Crouch1993}. Their methods require a large number
of exponentials (compared to the later developments) and introduce
specific order conditions for the coefficients.
Later, Munthe-Kaas~\cite{MuntheKaas1995,MuntheKaas1998,MUNTHEKAAS1999115}
introduced a class of
integrators that involve commutators and allow one to build a Lie group integrator
based on an arbitrary classical Runge-Kutta scheme. Then Celledoni, Marthinsen
and Owren~\cite{CELLEDONI2003341} developed another class of Lie group methods
that avoid commutators which results in a different structure of the coefficients
and the order conditions that complement the classical ones.
The complete theory of order conditions for commutator-free methods was worked
out by Owren in Ref.~\cite{Owren_2006}.

The main purpose of the present paper is to introduce a new class of
commutator-free Lie group methods that is naturally related to low-storage schemes
of Williamson~\cite{WILLIAMSON198048} and has different properties in terms
of exponentials reuse compared to the methods available in the literature.
In the way the exponentials are reused this class of methods is also
related to the multirate infinitesimal step (MIS) methods of
Knoth and collaborators~\cite{KNOTH1998327,Wensch2009}.
The first instance of a method that belongs to the new proposed class
in the literature is, to the best of our knowledge, the 3-stage third-order
coefficient scheme introduced by L\"{u}scher in Ref.~\cite{Luscher:2010iy}.

This paper is organized as follows. In Sec.~\ref{sec_rk} we review
classical Runge-Kutta integrators including low-storage schemes,
in Sec.~\ref{sec_lie} we review several types of Lie group integrators that
exist in the literature. In Sec.~\ref{sec_new} we propose a new class of
low-storage commutator-free Lie group integrators with reuse of exponentials
and prove that a 3-stage scheme in the new format is of order $p=3$ global accuracy.
We then formulate a conjecture about low-storage commutator-free Lie group
methods with more than three stages and of order higher than three.
In Sec.~\ref{sec_num} we provide numerical evidence in support of the
conjecture and conclude in Sec.~\ref{sec_concl}.

\section{Classical Runge-Kutta integrators}
\label{sec_rk}

We first review the well-known facts about explicit Runge-Kutta integrators
and low-storage schemes and introduce the notation that will be used in
the following.

\subsection{Definitions and notation}

Consider a first-order differential equation for a function $y(t)$
\begin{equation}
\label{eq_dydt}
\frac{dy}{dt}=f(t,y).
\end{equation}
A standard explicit $s$-stage Runge-Kutta (RK) scheme\footnote{It is
implied here and later on
that when the upper bound on the index in a sum
is smaller than the lower bound, the sum is set to 0 and if the same
conditions hold for a product, the product is set to 1, e.g. $\sum_{j=1}^0...=0$,
$\prod_{j=1}^0...=1$.}
for numerically integrating Eq.~(\ref{eq_dydt}) from time $t$ to $t+h$ 
is~\cite{ButcherBook,HairerBook1}
\begin{eqnarray}
y_i &=& y_t+h\sum_{j=1}^{i-1}a_{ij}k_j,\label{eq_yi}\\
k_i &=& f(t+hc_i,y_i),\label{eq_ki}\\
i&=&1,\dots,s,\\
y_{t+h}&=&y_t+h\sum_{i=1}^s b_ik_i.\label{eq_yth}
\end{eqnarray}
In the context of manifold integrators introduced later in Sec.~\ref{sec_lie}
we refer to this scheme as \textit{classical} RK method.
For an explicit method, $a_{ij}=0$ for $j\geqslant i$ and self-consistency conditions
require
\begin{equation}
\label{eq_ci}
c_i = \sum_{j=1}^{i-1} a_{ij}.
\end{equation}
The set of coefficients $a_{ij}$, $b_i$, $c_i$ can be conveniently displayed in a
Butcher table, for instance, for a 3-stage method:
\begin{equation}
\begin{array}
{c|lll}
c_2 & a_{21}\\
c_3 & a_{31} & a_{32} \\
\hline
& b_{1}   & b_{2}  & b_{3}  \\
\end{array}
\label{eq_abctable}
\end{equation}
where the first trivial entry $c_1=0$ is omitted.

Without loss of generality we focus on autonomous problems
\begin{equation}
\label{eq_dydt_a}
\frac{dy}{dt}=f(y).
\end{equation}
Extension to non-autonomous problems is straightforward and is not important
in the following.

By comparing the numerical solution (\ref{eq_yth}) with the Taylor expansion
of the exact solution $y(t+h)$ around $y(t)$ one obtains the constraints,
called \textit{the order conditions},
on the RK coefficients $a_{ij}$, $b_i$, $c_i$ so that the RK method provides
a certain order of accuracy. For example, the order conditions for a
$s$-stage RK method with global third-order accuracy are~\cite{ButcherBook,HairerBook1}
\begin{eqnarray}
\sum_i b_i &=& 1,\\
\sum_i b_ic_i &=& \frac{1}{2},\\
\sum_i b_ic_i^2 &=& \frac{1}{3},\\
\sum_{i,j} b_i a_{ij} c_j &=& \frac{1}{6}.
\end{eqnarray}
The minimum number of stages for a third-order RK method is three and the order
conditions then take the following form
\begin{eqnarray}
b_1 + b_2 + b_3 &=& 1,\label{eq_oc1}\\
b_2c_2 + b_3c_3 &=& \frac{1}{2},\\
b_2c_2^2 + b_3c_3^2 &=& \frac{1}{3},\\
b_3 a_{32} c_2 &=& \frac{1}{6}.\label{eq_oc4}
\end{eqnarray}
Given that there are six $a_{ij}$, $b_i$ coefficients (the coefficients $c_i$ follow
from (\ref{eq_ci})) and four constraints (\ref{eq_oc1})--(\ref{eq_oc4}), one expects
a two-parameter family of solutions. Due to the fact that $c_2$ and $c_3$ enter in
the constraints nonlinearly, it is customary to take $c_2$ and $c_3$ as free parameters and
in this case, there are three branches of solutions. Picking the most generic branch
$c_2\neq 0 \neq c_3\neq c_2\neq 2/3$~\cite{ButcherBook} one gets
\begin{eqnarray}
a_{32} &=& \frac{c_3(c_3-c_2)}{c_2(2-3c_2)},\label{eq_a32}\\
b_2 &=& \frac{3c_3-2}{6c_2(c_3-c_2)},\\
b_3 &=& \frac{2-3c_2}{6c_3(c_3-c_2)},\label{eq_b3}
\end{eqnarray}
and the other coefficients can be reconstructed trivially from the order conditions.

\subsection{Williamson low-storage schemes}
\label{sec_Wlow}

It was noted by Williamson in Ref.~\cite{WILLIAMSON198048} that the RK scheme 
(\ref{eq_yi})--(\ref{eq_yth}) with imposing additional constraints and a 
suitable choice of coefficients $A_i$, $B_i$ can be rewritten as ($A_1=0$)
\begin{eqnarray}
\Delta y_i &=& A_i\Delta y_{i-1}+hf(y_{i-1}),\label{eq_WRK1}\label{eq_Wdy}\\
y_i &=& y_{i-1} + B_i\Delta y_i,\\
i &=& 1,\dots,s,\\
y_0&\equiv& y_t,\\
y_{t+h}&\equiv&y_s.\label{eq_WRK2}
\end{eqnarray}
The utility of the scheme (\ref{eq_WRK1})--(\ref{eq_WRK2}) is that at a given
stage $i$ one only needs to keep the values of $y_i$ and $\Delta y_i$ and the previous
values can be overwritten. For a system of $N$ degrees of freedom only $2N$
storage locations are required, independently of the order and the number of stages
of the RK method. This particular two-register
low-storage scheme is referred to as $2N$-storage
scheme. The original RK coefficients are related to
$A_i$, $B_i$ as
\begin{eqnarray}
a_{ij} &=& \left\{
\begin{array}{ll}
A_{j+1}a_{i,j+1}+B_j, & j<i-1,\\
B_j, & j=i-1,\\
0, & \mbox{otherwise},
\end{array}
\right.\label{eq_Ai}\\
b_i &=&\left\{
\begin{array}{ll}
A_{i+1}b_{i+1}+B_i, & i<s,\\
B_i, & i=s.
\end{array}
\right.\label{eq_Bi}
\end{eqnarray}
The $2N$-storage schemes of Williamson~\cite{WILLIAMSON198048} have been modified
in various ways leading to the development of $2R$, $2S$, $3R$, etc.
schemes~\cite{KENNEDY2000177,KETCHESON20101763}
that differ in the number of registers (quantities stored at each
stage) and the constraints imposed on the coefficients $A_i$, $B_i$.
However, it is the $2N$-storage schemes that possess the properties that this
discussion builds upon later, so we consider only them here.

\begin{figure}[t]
	\centering
	\includegraphics[width=0.7\textwidth]{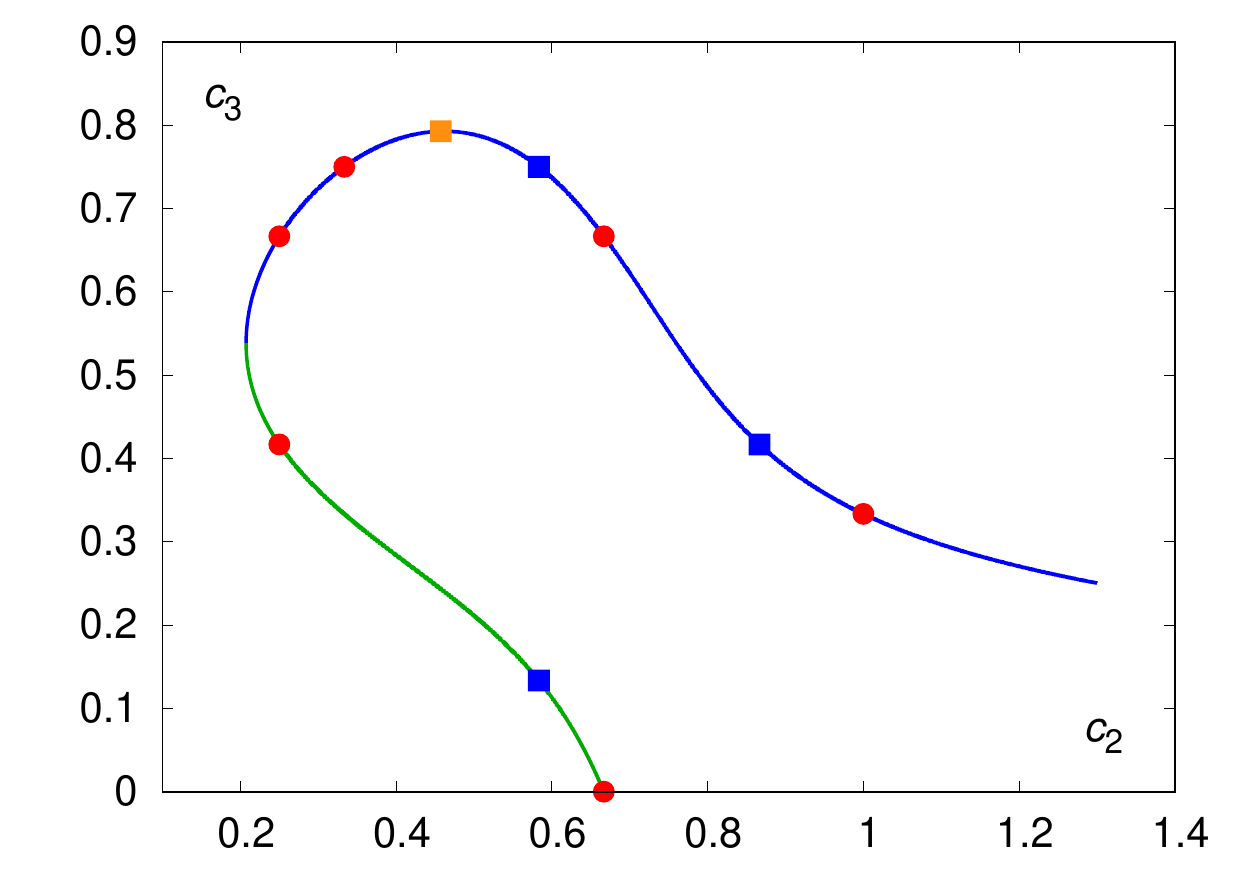}
	\caption{
		The two branches of solutions~\cite{WILLIAMSON198048}
		of Eq.~(\ref{eq_Wc3}) shown as blue and
		green lines.
		The red circles show the solutions with rational coefficients found
		in~\cite{WILLIAMSON198048} and the blue squares the solutions
		$(c_2=7/12,c_3=2/15)$, $(c_2=7/12,c_3=3/4)$ and $(c_2=13/15,c_3=5/12)$
		that were missed in~\cite{WILLIAMSON198048}.
		The orange square corresponds to the solution with minimal truncation
		error as defined in~\cite{Ralston1962} which is used later as
		a basis of the BWRRK33 commutator-free Lie group method,
		see Sec.~\ref{sec_num}.
		\label{fig_wrk_curves}
	}
\end{figure}

The $2N$-storage scheme introduces more constraints on the coefficients
$a_{ij}$, $b_i$. However, they may be implicit: once the classical coefficients
$a_{ij}$, $b_i$ are expressed in terms of the $2N$-storage coefficients $A_i$, $B_i$,
one needs to search for a solution that satisfies only the original
classical order conditions.
For low-order schemes it may be useful to find the additional constraints explicitly,
and we discuss a 3-stage third-order RK method in detail here to illustrate
this point. In this case the coefficients no longer form a two-parameter family.
The additional constraint can be imposed in different ways and a particular form
used in Ref.~\cite{WILLIAMSON198048} is
\begin{equation}
c_3^2(1-c_2)+c_3\left(c_2^2+\frac{1}{2}c_2-1\right)
+\left(\frac{1}{3}-\frac{1}{2}c_2\right)=0.\label{eq_Wc3}
\end{equation}
Choosing $c_2$ and then solving for $c_3$ from Eq.~(\ref{eq_Wc3}) allows one
to reconstruct the coefficients $a_{ij}$, $b_i$ from
(\ref{eq_a32})--(\ref{eq_b3}), (\ref{eq_oc1})--(\ref{eq_oc4}) and (\ref{eq_ci}).
Inverting the dependence (\ref{eq_Ai}), (\ref{eq_Bi}) produces the coefficients
$A_i$, $B_i$ of the $2N$-storage scheme (\ref{eq_Wdy})--(\ref{eq_WRK2}).

The two branches of solutions of $c_3$ as function of $c_2$
resulting from (\ref{eq_Wc3}) are shown
in Fig.~\ref{fig_wrk_curves}. There is a reflection symmetry 
with respect to the $c_2+c_3=1$
axis which is apparent after a change of variables $c_2=x+y$,
$c_3=1-x+y$.

\section{Brief review of Lie group integrators with examples at third order}
\label{sec_lie}

Let us now consider an equation of the form
\begin{equation}
\label{eq_dYAY}
\frac{dY}{dt}=A(Y)Y,
\end{equation}
where $Y$ is a vector or a matrix. We use capital letters to emphasize that
we now deal with objects that may not necessarily commute.
Again, for simplicity we consider autonomous problems and extension
to non-autonomous problems with $A(t,Y)$ is straightforward.
Although the primary focus in this article is equations on Lie groups where
$Y\in G$ and $A(Y)\in \mathfrak{g}$ where $G$ is a (matrix) Lie group and
$\mathfrak{g}$ its Lie algebra, the discussion applies to structure-preserving
integration of differential equations on manifolds
in general~\cite{Hairer2006}. 
The numerical examples in Sec.~\ref{sec_num} include
free rigid body rotation, where $Y$ is a three-dimensional vector of fixed
length and the manifold is $S^2$, integration of the gradient flow on $SU(3)$
where $Y$ is an $SU(3)$ matrix and the manifold is obviously the $SU(3)$ group,
van der Pol oscillator where $Y$ is a two-dimensional vector and the manifold
is $R^2$ and more.

In the next subsections we review several existing Lie group integrators
with examples at third order to define the building blocks necessary for the
discussion in Sec.~\ref{sec_new}.

\subsection{Crouch-Grossman methods}

An update from $Y_t$ to $Y_{t+h}$ in the form of a classical RK method
(\ref{eq_ki})--(\ref{eq_yth}) is possible, however, even if
$Y\in G$, the updated value of the form $Y+Const\cdot hA(Y)Y$ 
is no longer in the Lie group, in general. To maintain $Y$ on the manifold
one needs to construct an update of the form $\exp(Const\cdot hA(Y))Y$.
Then every stage of the RK-based method and the resulting $Y_{t+h}$ stays
on the original manifold.

Crouch and Grossman~\cite{Crouch1993} suggested an $s$-stage Lie group RK
method of the following form:
\begin{eqnarray}
Y_i &=& {\cal T}\left\{\prod_{j=1}^{i-1}\exp(h a_{ij}K_j)\right\}Y_t,\\
K_i &=& A(Y_i),\label{eq_CGKi}\\
i&=&1,\dots,s,\\
Y_{t+h}&=&{\cal T}\left\{\prod_{i=1}^{s}\exp(h b_{i}K_i)\right\}Y_t.\label{eq_CGYth}
\end{eqnarray}
Here ${\cal T}\prod$ represents a ``time-ordered''\footnote{by analogy with
quantum field theory} product with a convention that
 an element with \textit{smaller} value of the
index is always located \textit{to the right}. An explicit example 
clarifying this notation is given below in
Eq.~(\ref{eq_CG3Y1})--(\ref{eq_CG3Yth}).
In this case the number of $a_{ij}$, $b_i$ coefficients matches a classical
$s$-stage RK method and can be represented with a Butcher table.
The coefficients need to satisfy the classical order conditions and
some additional constraints that result from non-commutativity.
Ref.~\cite{Crouch1993} considered methods up to order three and found
that for a third-order method one has the following additional relation
\begin{equation}
\label{eq_CGoc}
\sum_i b_i^2 c_i + 2\sum_{i<j} b_i c_i b_j = \frac{1}{3}.
\end{equation}
The order conditions for higher-order Lie group methods of this type
were later derived in~\cite{OwrenMarthinsen1999}. It turns out that a
3-stage third-order RK Lie group method is possible
(since there are six coefficients and five constraints, giving a
one-parameter family)
and sets of coefficients
satisfying (\ref{eq_oc1})--(\ref{eq_oc4}) and (\ref{eq_CGoc}) were
given in Ref.~\cite{Crouch1993}. Let us for the later convenience
write the 3-stage third-order Crouch-Grossman RK method explicitly:
\begin{eqnarray}
Y_1 &=& Y_t,\label{eq_CG3Y1}\\
K_1 &=& A(Y_1),\\
Y_2 &=& \exp(h a_{21}K_1)Y_t,\\
K_2 &=& A(Y_2),\\
Y_3 &=& \exp(h a_{32}K_2)\exp(h a_{31}K_1)Y_t,\\
K_3 &=& A(Y_3),\\
Y_{t+h} &=& \exp(h b_{3}K_3)\exp(h b_{2}K_2)\exp(h b_{1}K_1)Y_t.\label{eq_CG3Yth}
\end{eqnarray}
From the computational perspective one can note the following. The method
requires three evaluations of the right hand side of the differential
equation, six exponentiations and storage for $K_i$ from all three stages,
to be applied at the last step of the algorithm (\ref{eq_CG3Yth}).

\subsection{Munthe-Kaas methods}

Another direction in constructing Lie group methods was taken by Munthe-Kaas
in Refs.~\cite{MuntheKaas1995,MuntheKaas1998,MUNTHEKAAS1999115}.
The most general approach worked out in
Ref.~\cite{MUNTHEKAAS1999115} represents the solution $Y(t)$ as $Y(t)=\exp(U(t))Y(0)$
and constructs an algorithm for solving the equation for $U(t)$
\begin{equation}
\frac{d U}{dt}=d\exp_U^{-1}(A(Y(t))),
\end{equation}
where the inverse derivative of the matrix exponential can be written
as an expansion
\begin{equation}
d\exp_U^{-1}=\sum_{k=0}^\infty\frac{B_k}{k!}{\rm ad}_U^k.
\end{equation}
$B_k$ are the Bernoulli numbers and the adjoint operator ${\rm ad}_U$ represents
a mapping ${\rm ad}_U(V)=[U,V]=UV - VU$. The $k$-th power of ${\rm ad}_U$ is understood
as an iterated application of this mapping:
\begin{eqnarray}
{\rm ad}_U^0(V)&=&V,\\
{\rm ad}_U^1(V)&=&[U,V],\\
{\rm ad}_U^k(V)&=&{\rm ad}_U({\rm ad}_U^{k-1}(V))=[U,[U,[\dots,[U,V]]]].\\
\end{eqnarray}
Let a truncated approximation of $d\exp_U^{-1}(V)$ be
\begin{equation}
\label{eq_dexpinv}
{\rm dexpinv}(U,V,p)=\sum_{k=0}^{p-1}\frac{B_k}{k!}{\rm ad}_U^k(V).
\end{equation}
Then using the notation introduced in earlier sections, a Lie group
$s$-stage order-$p$ RK method of Munthe-Kaas type has the following form:
\begin{eqnarray}
U_i &=& h\sum_{j=1}^{i-1}a_{ij}\tilde K_j,\label{eq_MKUi}\\
Y_i &=& \exp(U_i)Y_t,\\
K_i &=& A(Y_i),\\
\tilde K_i &=& {\rm dexpinv}(U_i,K_i,p),\label{eq_MKdexp}\\
i&=&1,\dots,s,\\
V &=& h\sum_{i=1}^s b_i\tilde K_i,\\
Y_{t+h}&=&\exp(V)Y_t.\label{eq_MKYth}
\end{eqnarray}
As shown in Ref.~\cite{MUNTHEKAAS1999115}, if the coefficients $a_{ij}$, $b_i$ correspond
to a classical RK method of order $p$ then the algorithm
(\ref{eq_MKUi})--(\ref{eq_MKYth}) with the truncation at $p-1$ in
(\ref{eq_dexpinv}) is a Lie group integrator of order at least $p$.
This procedure allows one to turn any classical $s$-stage order $p$
RK method into a Lie group
integrator with the same number of stages and the same order
at the expense of introducing commutators at every stage, Eq.~(\ref{eq_MKdexp}).
The number of commutators can be reduced as discussed in~\cite{MuntheKaasOwren1999},
and for later comparisons we write explicitly an earlier 
version~\cite{MuntheKaas1998}
of the 3-stage third-order Lie group RK method of Munthe-Kaas type that requires
only one commutator at the final stage:
\begin{eqnarray}
Y_1 &=& Y_t,\\
K_1 &=& A(Y_1),\\
Y_2 &=& \exp(h a_{21}K_1)Y_t,\\
K_2 &=& A(Y_2),\\
Y_3 &=& \exp(h (a_{32}K_2 + a_{31}K_1))Y_t,\\
K_3 &=& A(Y_3),\\
V &=& h\sum_{i=1}^3 b_{i}K_i,\\
\tilde V &=& V - \frac{h}{6}[K_1,V],\label{eq_MK3com}\\
Y_{t+h} &=& \exp(\tilde V)Y_t.
\end{eqnarray}
Apart from the exponential action and the commutator in (\ref{eq_MK3com})
this method resembles a classical RK method in that respect that one
adds $K_i$ in a similar fashion as in a classical method and then exponentiates
the result to produce $Y_i$ for the next stage. Here one needs three evaluations
of the right hand side, three exponentiations and storage of $K_i$ from all 
three stages.

\subsection{Celledoni-Marthinsen-Owren methods}

Celledoni, Marthinsen and Owren in Ref.~\cite{CELLEDONI2003341} considered an approach
that generalizes Crouch-Grossman methods with the goal to avoid computation
of commutators. Their idea is to introduce more than one exponential per
stage of a RK method but allow for linear combinations of $K_i$ in the
exponentials. An $s$-stage RK Lie group integrator can
be written in the following form:
\begin{eqnarray}
Y_i &=& {\cal T}\left\{\prod_{l=1}^{L_i}
\exp\left(h \sum_{j=1}^{J_{il}} \alpha_{l;ij}K_j\right)\right\}Y_t,\label{eq_CMOYi}\\
K_i &=& A(Y_i),\\
i&=&1,\dots,s,\\
Y_{t+h} &=& {\cal T}\left\{\prod_{l=1}^{L}
\exp\left(h \sum_{i=1}^{I_l} \beta_{l;i}K_i\right)\right\}Y_t.\label{eq_CMOYth}
\end{eqnarray}
The notation here is similar but not the same as in Ref.~\cite{CELLEDONI2003341}
to be more in line with the methods introduced earlier. Here
$L_i$ is the number of exponentials used at stage $i$, $J_{il}$ is the upper
bound on summation inside the $l$-th exponential at $i$-th stage, and,
similarly, $L$ is the number of exponentials at the final stage and $I_l$
is the upper bound on summation inside the $l$-th exponential at the final
stage. By introducing more parameters one has more room to satisfy the
additional order conditions arising from non-commutativity at the expense
of introducing more exponentials at each stage. The new coefficients
$\alpha_{l;ij}$, $\beta_{l;i}$ are related to the coefficients of a classical
RK method as~\cite{CELLEDONI2003341}
\begin{eqnarray}
\sum_{l=1}^{L_i}\alpha_{l;ij}=a_{ij},\label{eq_aCMO}\\
\sum_{l=1}^{L}\beta_{l;i}=b_{i}.\label{eq_bCMO}
\end{eqnarray}
The Crouch-Grossman method, Eqs.~(\ref{eq_CG3Y1})--(\ref{eq_CG3Yth})
is a subclass of the 
Celledoni-Marthinsen-Owren methods where $\alpha_{l;ij}=a_{ij}\delta_{lj}$
(and automatically $L_i=i-1$ for explicit methods).

Ref.~\cite{CELLEDONI2003341} proceeded in a way that minimizes the number
of exponentials and constructed schemes of third and fourth order
that have the minimal number of exponentials 
and also reuse the exponentials at next stages. Here for
comparison we write explicitly one of the solutions found in~\cite{CELLEDONI2003341}:
\begin{eqnarray}
Y_1&=&Y_t,\\
K_1&=&A(Y_1),\\
Y_2&=&\exp(h\alpha_{1;21}K_1)Y_t,\\
K_2&=&A(Y_2),\\
Y_3&=&\exp(h(\alpha_{1;32}K_2+\alpha_{1;31}K_1))Y_t,\\
K_3&=&A(Y_3),\\
Y_{t+h} &=& \exp(h(\beta_{2;3}K_3+\beta_{2;2}K_2+\beta_{2;1}K_1))
\exp(h\beta_{1;1}K_1)Y_t.
\end{eqnarray}
Requiring that $\beta_{1;1}=\alpha_{1;21}$ allows one to reuse $Y_2$ and
calculate only one exponential at the last stage. This requirement also
fixes these two coefficients to be equal to $1/3$ and the other coefficients
then form a one-parameter family of solutions and their explicit
form is given in Ref.~\cite{CELLEDONI2003341}. Another branch of solutions
reuses $Y_3$ and results in a method with the same computational requirements:
three right hand side evaluations, three exponentiations and storage of $K_i$
from all stages and $Y_2$ or $Y_3$.

As one can see, at third order the Munthe-Kaas and Celledoni-Marthinsen-Owren
methods have similar computational requirements, however, at fourth and higher order the
situation is different: while Munthe-Kaas method can be constructed with the same
number of exponentials as the number of stages in a classical RK method
(with one exponential per stage), the Celledoni-Marthinsen-Owren
methods require more exponentials (\textit{e.g.}, at least, five at fourth order).

\section{A new class of commutator-free Lie group integrators}
\label{sec_new}

\subsection{Construction of the integrator}
\label{sec_construct}

Here we construct a new class of Lie group integrators. It can be considered
as a subclass of Celledoni-Marthinsen-Owren methods,
however, the construction proceeds differently
and results in a family of solutions different from Ref.~\cite{CELLEDONI2003341}.
In particular this new scheme has different properties in terms of storage and
exponentials reuse.
At the same time, this new class can be considered as
a subclass of the multirate infinitesimal step (MIS) methods used as exponential
integrators~\cite{Wensch2009}, however, again, the family of solutions
proposed here is different from the ones present in the literature.

Let us first construct a 3-stage third-order Lie group integrator and then
comment on generalization of this scheme. Let us take the structure introduced
in Eq.~(\ref{eq_CMOYi})--(\ref{eq_CMOYth}) and add the following requirements:
\newcounter{store_list}
\begin{enumerate}
	\item $L_i=i-1$ -- as in the Crouch-Grossman method,
	stage $i$ has exactly $i-1$ exponentials.
	\item $J_{il}=l$ -- the number of terms within each exponential is equal
	to the index of that exponential in the sequence. With the time-ordering
	convention this means that the rightmost exponential has one term,
	the one to the left of it -- two terms, and so on.
	\item $L=s$ -- at the final stage there is the maximum number, $s$ exponentials.
	\item $I_l=l$ -- the convention on the number of terms inside exponentials
	at the final stage is the same as in the previous stages.
	\setcounter{store_list}{\value{enumi}}
\end{enumerate}

Explicitly, a 3-stage algorithm (its order is not yet determined) is
\begin{eqnarray}
Y_1&=&Y_t,\label{eq_exp3_y1}\\
K_1&=&A(Y_1),\\
Y_2&=&\exp(h\alpha_{1;21}K_1)Y_t,\\
K_2&=&A(Y_2),\\
Y_3&=&\exp(h(\alpha_{2;32}K_2+\alpha_{2;31}K_1))
\exp(h\alpha_{1;31}K_1)Y_t,\\
K_3&=&A(Y_3),\\
Y_{t+h} &=& 
\exp(h(\beta_{3;3}K_3+\beta_{3;2}K_2+\beta_{3;1}K_1))\\
&\times&
\exp(h(\beta_{2;2}K_2+\beta_{2;1}K_1))
\exp(h\beta_{1;1}K_1)Y_t.\label{eq_exp3_yth}
\end{eqnarray}
This algorithm has six exponentials as the Crouch-Grossman method and also
requires evaluating linear combinations of $K_i$ as in a classical RK
method. There are 10 coefficients $\alpha_{l;ij}$, $\beta_{l;i}$ that are related
to the classical RK coefficients via (\ref{eq_aCMO}) and (\ref{eq_bCMO}) and
are subject to the four classical order conditions
(\ref{eq_oc1})--(\ref{eq_oc4})
and possibly other constraints arising from noncommutativity.

At first sight,
there is nothing beneficial in this scheme as it requires more work
than any other Lie group method introduced previously and therefore we apply
another constraint:
\begin{enumerate}
	\setcounter{enumi}{\value{store_list}}
	\item The coefficients in the exponentials with the same number of terms
	are the same at all stages.
\end{enumerate}
This means $\beta_{1;1}=\alpha_{1;31}=\alpha_{1;21}$, $\beta_{2;1}=\alpha_{2;31}$
and $\beta_{2;2}=\alpha_{2;32}$.
This requirement allows one to reuse previous $Y_i$ at every stage 
and the scheme can be rewritten as
\begin{eqnarray}
Y_1&=&Y_t,\label{eq_myY1}\\
K_1&=&A(Y_1),\\
Y_2&=&\exp(h\alpha_{1;21}K_1)Y_1,\\
K_2&=&A(Y_2),\\
Y_3&=&\exp(h(\alpha_{2;32}K_2+\alpha_{2;31}K_1))Y_2,\\
K_3&=&A(Y_3),\\
Y_{t+h} &=& 
\exp(h(\beta_{3;3}K_3+\beta_{3;2}K_2+\beta_{3;1}K_1))Y_3.\label{eq_myYth}
\end{eqnarray}
Now there are only three exponentials, the method reuses values of $Y_i$ from
each previous stage and if the coefficients can be tuned that the scheme results in
a third-order method, it can be on par with the methods of Sec.~\ref{sec_lie}.
There are now six independent coefficients, as in the classical 3-stage RK
method and they are related to the coefficients of the classical method in a simple
way:
\begin{eqnarray}
\alpha_{1;21} &=& a_{21},\\
\alpha_{2;31} &=& a_{31}-a_{21},\\
\alpha_{2;32} &=& a_{32},\\
\beta_{3;1} &=& b_1-a_{31},\\
\beta_{3;2} &=& b_2-a_{32},\\
\beta_{3;3} &=& b_3.
\end{eqnarray}

Note that although a 3-stage third-order method is considered here as an example,
the construction (\ref{eq_myY1})--(\ref{eq_myYth}) is applicable in general.
Eqs.~(\ref{eq_CMOYi})--(\ref{eq_CMOYth}) with the five requirements listed
above essentially mean that 
in this format for a $s$-stage method
each stage $i$ has only one exponential that contains
a sum of all $K_i$ accumulated up to that stage that multiplies $Y_{i-1}$ from the
previous stage:
\begin{eqnarray}
Y_0&\equiv&Y_t,\label{eq_myY0gen}\\
Y_i &=& 
\exp\left(h \sum_{j=1}^{i-1} \alpha_{i-1;ij}K_j\right)Y_{i-1},\\
K_i &=& A(Y_i),\\
i&=&1,\dots,s,\\
Y_{t+h} &=& 
\exp\left(h \sum_{i=1}^{s} \beta_{s;i}K_i\right)Y_s.\label{eq_myYthgen}
\end{eqnarray}
However, as will be shown immediately below, a more compact
format may be possible.

\subsection{Order conditions for the new three-stage third-order Lie group method}
\label{sec_myoc}

By Taylor expanding the scheme (\ref{eq_myY1})--(\ref{eq_myYth}) and
comparing with the expansion of the exact solution one finds that the
additional order conditions for this scheme to be globally of third order can be
written as
\begin{equation}
b_2c_2^2+b_3c_3^2+(b_2c_2+b_3c_3)(b_1+b_2+b_3+c_3)+a_{32}c_2(c_2-b_1-b_2-b_3)=1,
\end{equation}
\begin{equation}
(b_2c_2+b_3c_3)(b_1+b_2+b_3-c_3)+a_{32}c_2(b_1+b_2+b_3-c_2)=\frac{1}{3}.
\end{equation}
With the use of the classical order conditions (\ref{eq_oc1})--(\ref{eq_oc4})
one finds that these two conditions are not independent and result in a
single condition:
\begin{equation}
\label{eq_ocmy}
a_{32}c_2(1-c_2)=\frac{1}{6}(3c_3-1).
\end{equation}
We can now multiply both sides by $b_3\neq0$, use Eqs.~(\ref{eq_oc4}) 
and (\ref{eq_b3}) and rewrite (\ref{eq_ocmy}) as a relation between $c_2$ and $c_3$:
\begin{equation}
c_3^2(1-c_2)+c_3\left(c_2^2+\frac{1}{2}c_2-1\right)
+\left(\frac{1}{3}-\frac{1}{2}c_2\right)=0.\label{eq_Wc3_again}
\end{equation}
Eq.~(\ref{eq_Wc3_again}) is exactly the same (!) as the relation (\ref{eq_Wc3})
for the $2N$-storage scheme discussed in Sec.~\ref{sec_Wlow}.

Let us summarize what has been achieved so far.
We proposed a new format for a 3-stage
commutator-free Lie group RK method (\ref{eq_exp3_y1})--(\ref{eq_exp3_yth}) and,
by requiring that it is of third order global accuracy, found that the order conditions
on the classical RK coefficients of this method are the same as on the 
3-stage third-order $2N$-storage scheme. Note, that the $2N$-storage
schemes~\cite{WILLIAMSON198048} were not intended as Lie group integrators and
were designed as classical RK methods. The other way around, this means that
although it was not imposed in (\ref{eq_myY1})--(\ref{eq_myYth}), 
the relations between the coefficients are such that one does not need
to store $K_i$ from all stages and this scheme
can be rewritten in a $2N$-storage format 
by analogy with
(\ref{eq_Wdy})--(\ref{eq_WRK2})
as\footnote{The clash of notation here
is unfortunate, but it is customary in the literature on low-storage schemes
to use $A_i$ for the coefficients, while it is customary in the literature 
on Lie group methods to use $A(Y)$ on the right hand side 
of Eq.~(\ref{eq_dYAY}).} ($A_1=0$)
\begin{eqnarray}
\Delta Y_i &=& A_i\Delta Y_{i-1}+hA(Y_{i-1}),\label{eq_mydy}\\
Y_i &=& \exp(B_i\Delta Y_i)Y_{i-1},\\
i &=& 1,\dots,s,\\
Y_0&\equiv& Y_t,\\
Y_{t+h}&\equiv&Y_s,\label{eq_myY}
\end{eqnarray}
where the coefficients $A_i$, $B_i$ are related to $a_{ij}$, $b_i$
in Eqs.~(\ref{eq_Ai}) and (\ref{eq_Bi}).

\subsection{Conjecture about higher order 
	commutator-free Lie group integrators}
\label{sec_newconj}

The first main result of this paper derived in 
the previous section can be summarized as follows:
\textit{The family of classical $2N$-storage 3-stage third-order RK
schemes are also automatically third-order commutator-free Lie group
integrators.}

In fact, the third-order scheme, presented without derivation in the
Appendix of Ref.~\cite{Luscher:2010iy} and used as a Lie group method
for integration of $SU(3)$ gradient flow, belongs to the class 
of integrators proposed in Sec.~\ref{sec_construct} with a specific
choice of coefficients from the one-parameter family of
Eq.~(\ref{eq_Wc3_again}). This is discussed in more detail in the third
numerical example in Sec.~\ref{sec_ex3}.

A natural question is: Are the $2N$-storage schemes at third order with
more than three stages and at orders four and higher also 
commutator-free Lie group methods of the same order? While there is no analytic proof
immediately available, the numerical evidence that is examined
in Sec.~\ref{sec_num} suggests that the answer to this question may be positive.

The second main result of this paper is the following conjecture:
\textit{$2N$-storage $s$-stage classical RK schemes of order $p$ are also automatically
commu\-tator-free Lie group methods
of the format proposed in Eqs.~(\ref{eq_mydy})--(\ref{eq_myY})
of order $p$ for orders $p=3$, $4$, $5$ and
possibly higher.}

\section{Numerical experiments}
\label{sec_num}

We now consider a few $2N$-storage classical RK schemes and apply them in
several examples to provide support for the conjecture stated in
Sec.~\ref{sec_newconj}.

First, we consider a 3-stage third-order family of $2N$-storage schemes for which
it is proven in Sec.~\ref{sec_myoc} that they are Lie group integrators
of order $p=3$. We would like to choose a set of coefficients for which the truncation
error is minimal in the sense of Ref.~\cite{Ralston1962}. We need to
emphasize the difference with the classical RK case. The 3-stage third-order 
classical RK scheme that has minimal truncation error found by
Ralston~\cite{Ralston1962} is not
a $2N$-storage scheme and is not of third order if used as a Lie group integrator
as defined in (\ref{eq_myY1})--(\ref{eq_myYth}). However, one can follow the error
minimization criteria of~\cite{Ralston1962} with the additional constraint
between $c_2$ and $c_3$, Eq.~(\ref{eq_Wc3_again}). The resulting set of
classical RK coefficients
was found by Williamson~\cite{WILLIAMSON198048} and they turn out to be not
rational. Unfortunately, Ref.~\cite{WILLIAMSON198048} provided the coefficients
with only five digits of accuracy which is not sufficient for the tests in this
section. Therefore we improve on this by following the minimization procedure
of~\cite{Ralston1962} with the constraint~(\ref{eq_Wc3_again}) and the resulting
set of coefficients is
\begin{eqnarray}
a_{21} &=& \phantom{-}0.45737999756938819,\label{eq_opt_a21}\\
a_{31} &=& -0.13267640849031470,\\
a_{32} &=& \phantom{-}0.92529641092092174,\\
b_1  &=& \phantom{-}0.19546562910003523,\\
b_2  &=& \phantom{-}0.41072077622489378,\\
b_3  &=& \phantom{-}0.39381359467507099.\label{eq_opt_b3}
\end{eqnarray}
By using (\ref{eq_Ai}), (\ref{eq_Bi}) they are transformed into the
$2N$-storage format and used in the commutator-free Lie group method,
Eqs.~(\ref{eq_mydy})--(\ref{eq_myY}). We use this scheme,
called BWRRK33, in the examples below,
and we also tested it
with the nine sets of rational coefficients shown in Fig.~\ref{fig_wrk_curves},
six of which were found in~\cite{WILLIAMSON198048}.

Next, we also consider thirteen $2N$-storage schemes available in the literature.
The nomenclature used here is the following. Letters ``RK'' in the middle
indicate that this is a classical RK method, the letters in front abbreviate
the names of the authors or the name given to the scheme in the original
article, the last digit is the order of the method, the digits in front of it
represent the number of stages and the additional letters after ``RK'' possible
notation from the original article to distinguish integrators with different
properties. The list of all fourteen $2N$-storage schemes tested in the examples is
given in Table~\ref{tab_ints}.

\begin{table}
\centering
\caption{
Method 1, BWRRK33 is proven to be a third-order commutator-free Lie group
method in Sec.~\ref{sec_myoc}. Methods 2--14 are classical 
$2N$-storage RK methods available in the literature used to test the conjecture
formulated in Sec.~\ref{sec_newconj}. We note that for the BWRRK33 method
we also ran numerical tests
with the rational coefficients shown in Fig.~\ref{fig_wrk_curves} and for the
CKRK54 method we ran tests with all four sets of coefficients found in
Ref.~\cite{CK1994}. For presenting results we, however, chose only one, recommended set
of coefficients from~\cite{CK1994}. Similarly, there are four sets of coefficients
for the BPRKO73 method in Ref.~\cite{BERNARDINI20094182} which we tested, but for
presenting the results we chose the set called 
ORK37-6 in Ref.~\cite{BERNARDINI20094182}.
Note also that BBBRKNL64 is called RK46-NL in Ref.~\cite{BERLAND20061459}.
\label{tab_ints}
}
\begin{tabular}{rlccc}
 & Name & Stages & Order & Reference \\
\hline
1  & BWRRK33   &  3 & 3 & Here, \cite{WILLIAMSON198048}, \cite{Ralston1962}\\
2  & BPRKO73   &  7 & 3 & \cite{BERNARDINI20094182} \\
3  & TSRKC73   &  7 & 3 & \cite{TOULORGE20122067} \\
4  & CKRK54    &  5 & 4 & \cite{CK1994} \\
5  & SHRK64    &  6 & 4 & \cite{STANESCU1998674} \\
6  & BBBRKNL64 &  6 & 4 & \cite{BERLAND20061459} \\
7  & HALERK64  &  6 & 4 & \cite{ALLAMPALLI20093837} \\
8  & HALERK74  &  7 & 4 & \cite{ALLAMPALLI20093837} \\
9  & TSRKC84   &  8 & 4 & \cite{TOULORGE20122067} \\
10 & TSRKF84   &  8 & 4 & \cite{TOULORGE20122067} \\
11 & NDBRK124  & 12 & 4 & \cite{NIEGEMANN2012364} \\
12 & NDBRK134  & 13 & 4 & \cite{NIEGEMANN2012364} \\
13 & NDBRK144  & 14 & 4 & \cite{NIEGEMANN2012364} \\
14 & YRK135    & 13 & 5 & \cite{Yan2017}
\end{tabular}
\end{table}

It is important to stress that while we proved
that the BWRRK33 scheme\footnote{and all 3-stage third-order
explicit RK schemes
satisfying the constraint (\ref{eq_Wc3_again})}
is a commutator-free Lie group integrator with 
global accuracy of order $p=3$,  Sec.~\ref{sec_myoc}, none of the other schemes in Table~\ref{tab_ints}
were originally designed as Lie group integrators. They are $2N$-storage
classical RK methods which were designed to have specific properties such as
increased stability regions, low dissipation, etc. If one attempts to use
an order $p$ arbitrary classical RK scheme whose coefficients satisfy only the
classical RK constraints as a Lie group integrator, the order of accuracy
is less than $p$, and as numerical experiments show, typically second order
at best. Thus, the thirteen $2N$-storage schemes (in addition to BWRRK33)
collected in Table~\ref{tab_ints}
are a perfect test of the conjecture in Sec.~\ref{sec_newconj} if they 
maintain the same order of accuracy when used as commutator-free Lie group
methods in the sense of Eqs.~(\ref{eq_mydy})--(\ref{eq_myY}).

Also, while it is not yet proven (or refuted) that any $2N$-storage method
of order $p$
is also a commutator-free Lie group integrator of order $p$,
for a given set of numerical values of the coefficients
the order conditions can be 
algorithmically checked by using B-series~\cite{Knoth2020git}.
All the methods of Table~\ref{tab_ints} were independently checked by Knoth
with the software available at~\cite{Knoth2020git} and they indeed fulfill
the order conditions corresponding to the order shown in the table when
used as Lie group integrators.

For numerical experiments in our code
we also implemented integrators of other types, 
namely, Crouch-Grossman
method of order $p=3$~\cite{Crouch1993}, Munthe-Kaas methods of order
$p=3$, $4$, $5$, $8$~\cite{MUNTHEKAAS1999115}, Celledoni-Marthinsen-Owren
methods~\cite{CELLEDONI2003341}
of order $p=3$ and $4$. In the tests below as a reference we use
the following Lie group integrators of Munthe-Kaas (RKMK) type:
\begin{itemize}
	\item 3-stage third-order with Ralston coefficients,
	\item 4-stage fourth-order with Ralston coefficients,
	\item 6-stage fifth-order with Butcher coefficients.
\end{itemize}

\subsection{Example 1: Free rigid body rotation}
\label{sec_ex1}

\begin{figure}[t]
	\centering
	\includegraphics[width=0.7\textwidth]{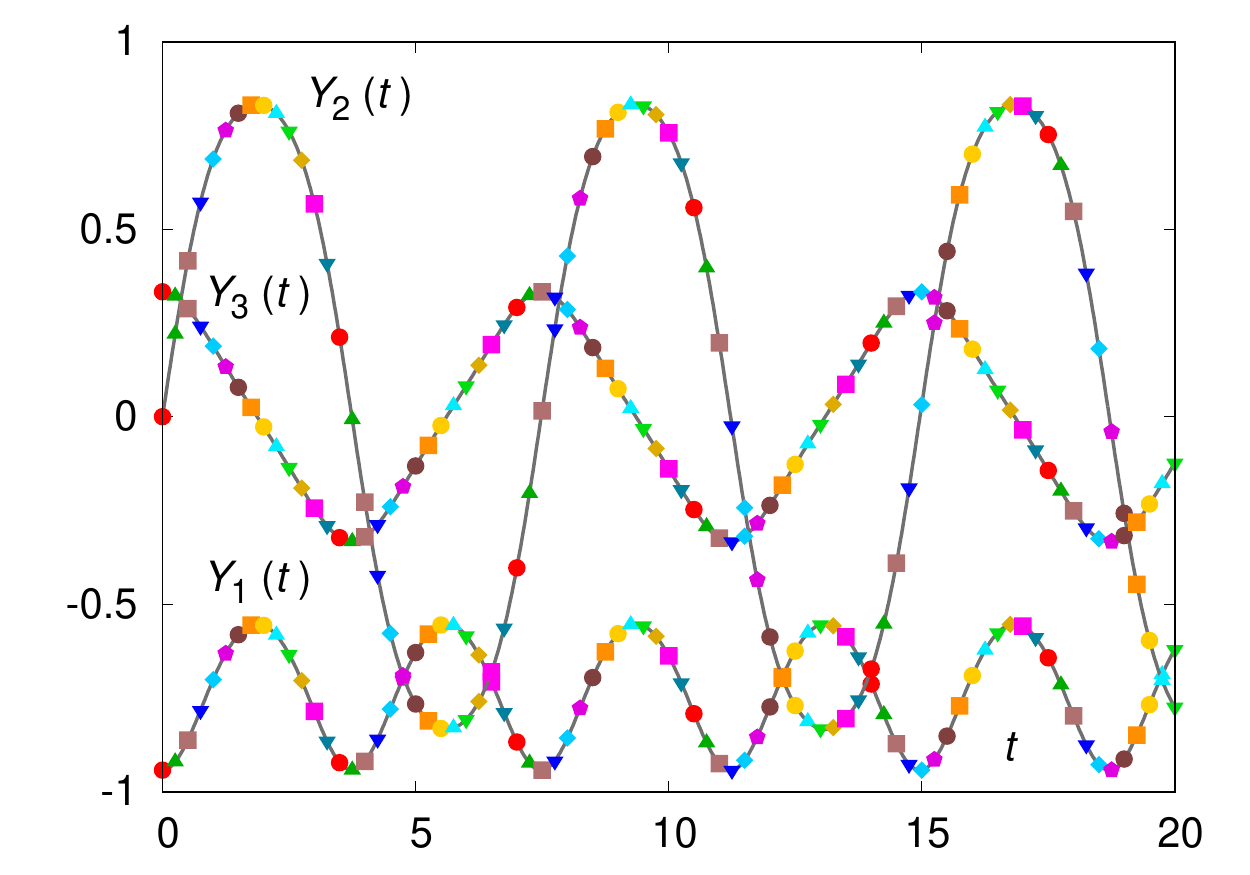}
	\caption{
		Comparison of the exact solution for the three components
		of the angular momentum (see text), shown as the gray lines,
		with the results produced by
		the fourteen integrators listed in Table~\ref{tab_ints}.
		If plotted as lines, all results are indistinguishable from the
		exact solution. Therefore we plot the results from different
		integrators as symbols of different shape and color skipping 140 steps in
		the sequence and starting to plot 
		the first integrator at a shift of 0 steps, 
		second at a shift of 10 steps and so on.
		\label{fig_ex1sol}
	}
\end{figure}

\begin{figure}[t]
	\centering
	\includegraphics[width=0.7\textwidth]{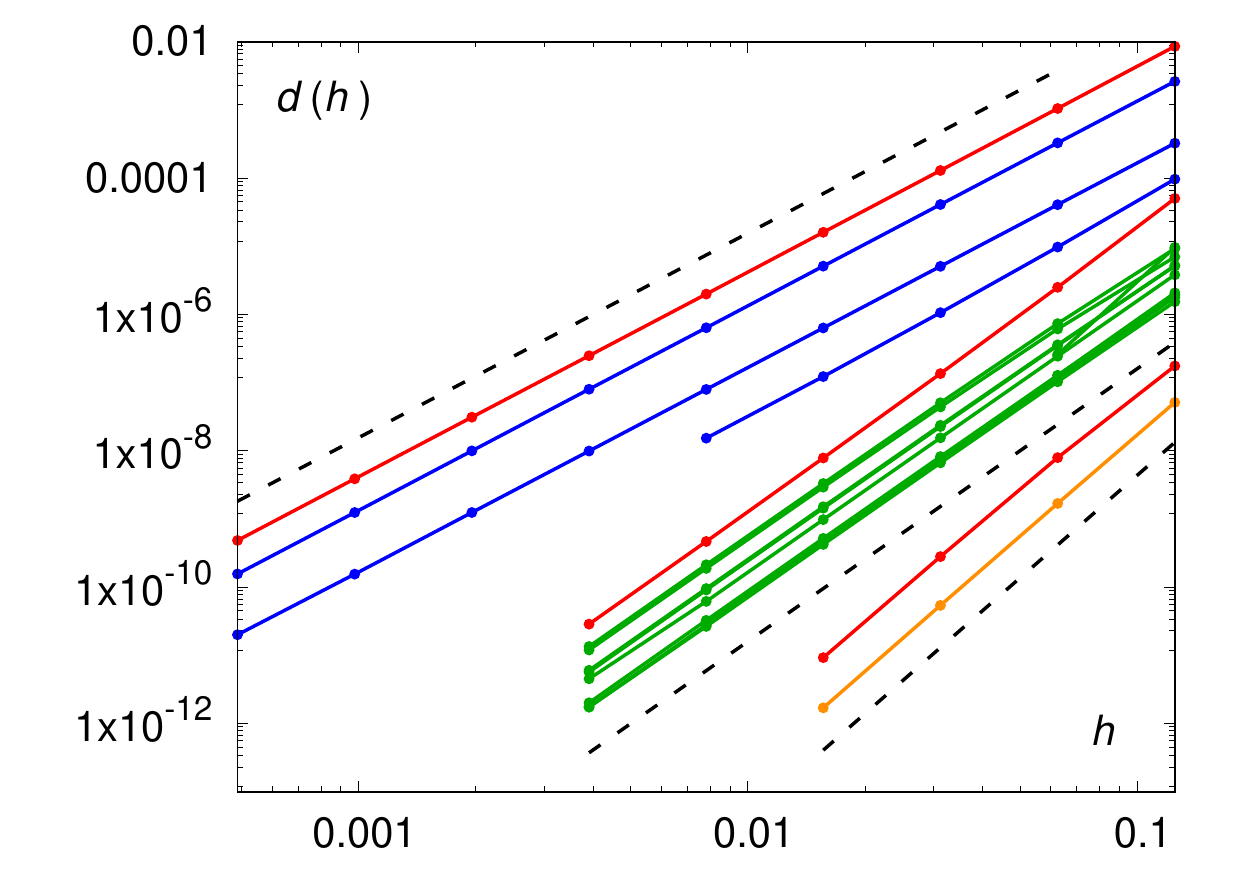}
	\caption{
		Distance from the reference (exact) solution $d(h)$ for various
		integrators as function of step size $h$ shown in a log-log plot
		for the rigid body problem, Eq.~(\ref{eq_rigid}).
		The red lines represent
		the three integrators of Munthe-Kaas type of order $p=3$, $4$ and $5$.
		The three blue lines represent the three integrators of order $p=3$ from
		Table~\ref{tab_ints}, the ten green lines the ten integrators of
		order $p=4$ and the orange line one integrator of order $p=5$ from the
		same table.
		For the $p=4$ integrators the minimum step size shown is $1/256$ and
		for $p=5$ $1/64$ since when $d(h)$ becomes comparable to $10^{-13}$
		the roundoff errors prevent correct scaling behavior. The black dashed
		lines are shown to guide the eye and represent from top to bottom
		$h^3$, $h^4$ and $h^5$, respectively.
		\label{fig_ex1scale}
	}
\end{figure}

As a first numerical example we consider
rotation of a free rigid body with the center of mass
fixed at the origin. This example was used in Ref.~\cite{CELLEDONI2003341,Curry2019}.
In this case, $Y$ is a three-dimensional vector of angular momentum
and the Euler equation is
\begin{equation}
\frac{dY}{dt}=Y\times I^{-1}Y,
\end{equation}
where $I$ is the inertia tensor. By using the hat map $\hat{\,}{\,}$: $R^3\to\mathfrak{so}(3)$ defined as
\begin{equation}
V = \left(
\begin{array}{c}
v_1\\
v_2\\
v_3\\
\end{array}
\right)\to
\hat V = \left(
\begin{array}{rrr}
0    & -v_3 & v_2 \\
v_3  & 0    & -v_1 \\
-v_2 & v_1  & 0
\end{array}
\right)
\end{equation}
the Euler equation can be rewritten in the form (\ref{eq_dYAY})
\begin{equation}
\label{eq_rigid}
\frac{dY}{dt}=-\widehat{I^{-1}Y}Y
\end{equation}
with $A(Y)\equiv-\widehat{I^{-1}Y}$. For the tensor of inertia we take the
same value $I=diag(7/8,5/8,1/4)$ as in~\cite{CELLEDONI2003341} but choose a
different initial condition $Y(0)=(-\sqrt{8}/3,0,1/3)$. Such an initial condition matches
the simplifying assumptions of~\cite{Marsden1999} where the exact solution
is given in terms of Jacobi's elliptic functions.
First, to check the implementation of the integrators in our code, we compare
the trajectory integrated from $t=0$ to $t=20$ with the time step $h=0.025$ with
all fourteen integrators of Table~\ref{tab_ints} with the exact solution.
The result is shown in Fig.~\ref{fig_ex1sol}.

\begin{figure}[t]
	\centering
	\includegraphics[width=0.7\textwidth]{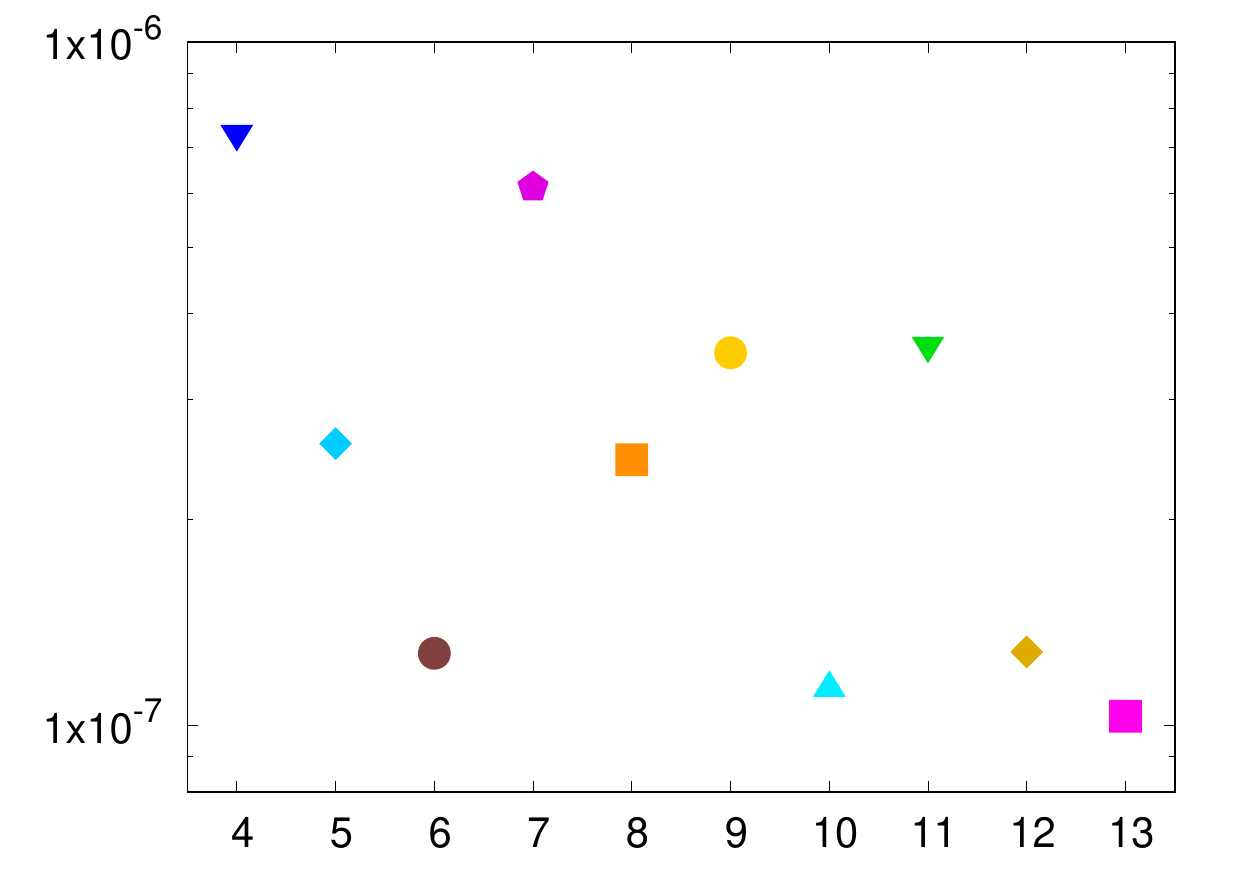}
	\caption{
		Distance from the reference solution $d(h)$ at $h=1/16$ 
		for all fourth-order integrators from Table~\ref{tab_ints}
		for the rigid body problem, Eq.~(\ref{eq_rigid}).
		The symbols and colors are the same as in Fig.~\ref{fig_ex1sol}
		and the numbers on the horizontal axis correspond to the
		numbering in Table~\ref{tab_ints}.
		\label{fig_ex1diffh16}
	}
\end{figure}

Next, to study the order of the methods we integrate the equation of motion
from $t=0$ to $t=3$ by using the step size $h=1/2^n$ where $n=3,\dots,11$.
Let $Y(t=3,h)$ be the solution evaluated at a particular step size $h$
and $Y_{ref}(t=3)$ the exact solution~\cite{Marsden1999}.
We define a distance metric as
\begin{equation}
d(h)=|Y(t=3,h)-Y_{ref}(t=3)|,
\end{equation}
where $|\dots|$ is the usual Euclidean vector norm. If an integrator has 
the global order of accuracy $p$ then one expects $d(h)\sim h^p$.
The results for $d(h)$ for the fourteen integrators of Table~\ref{tab_ints}
and the three reference integrators are shown in
Fig.~\ref{fig_ex1scale}. 
We note that the BPRKO73, SHRK64 and BBBRKNL64 integrators are somewhat
problematic since their coefficients are given with less than full double
precision accuracy. Their scaling breaks down when $d(h)$ reaches
about $10^{-8}$, $10^{-7}$ and $10^{-11}$, respectively.
Therefore we plot $d(h)$ approximately
down to those limits for those integrators.

As one can observe from Fig.~\ref{fig_ex1scale}, the $2N$-storage classical
RK schemes provide the same global order of accuracy when used as manifold
integrators of a new format defined in Eqs.~(\ref{eq_mydy})--(\ref{eq_myY}),
supporting the conjecture stated in Sec.~\ref{sec_newconj}.

Some of the RK methods of fourth order shown in Fig.~\ref{fig_ex1scale}
as green lines have comparable global errors and their lines are hard to
distinguish on the scale of the plot.
In Fig.~\ref{fig_ex1diffh16} we show the distance from the exact
solution $d(h)$ at $h=1/16$ for each method of
order 4 from Table~\ref{tab_ints}. For instance, methods 4 and 7,
5 and 8,
9 and 11, and 6, 10, 12 and 13 produce similar $d(h)$ in this example.
Similar features are observed in the other examples considered below
but which particular methods produce close results depends on the
differential equation.

A simple Matlab script illustrating the usage of BWRRK33, 
TSRKF84 and YRK135 as Lie group integrators for this example
is given in \ref{app_matlab}.

\subsection{Example 2: $SO(5)$ from Ref.~\cite{MUNTHEKAAS1999115}}
\label{sec_ex2}

\begin{figure}[t]
	\centering
	\includegraphics[width=0.7\textwidth]{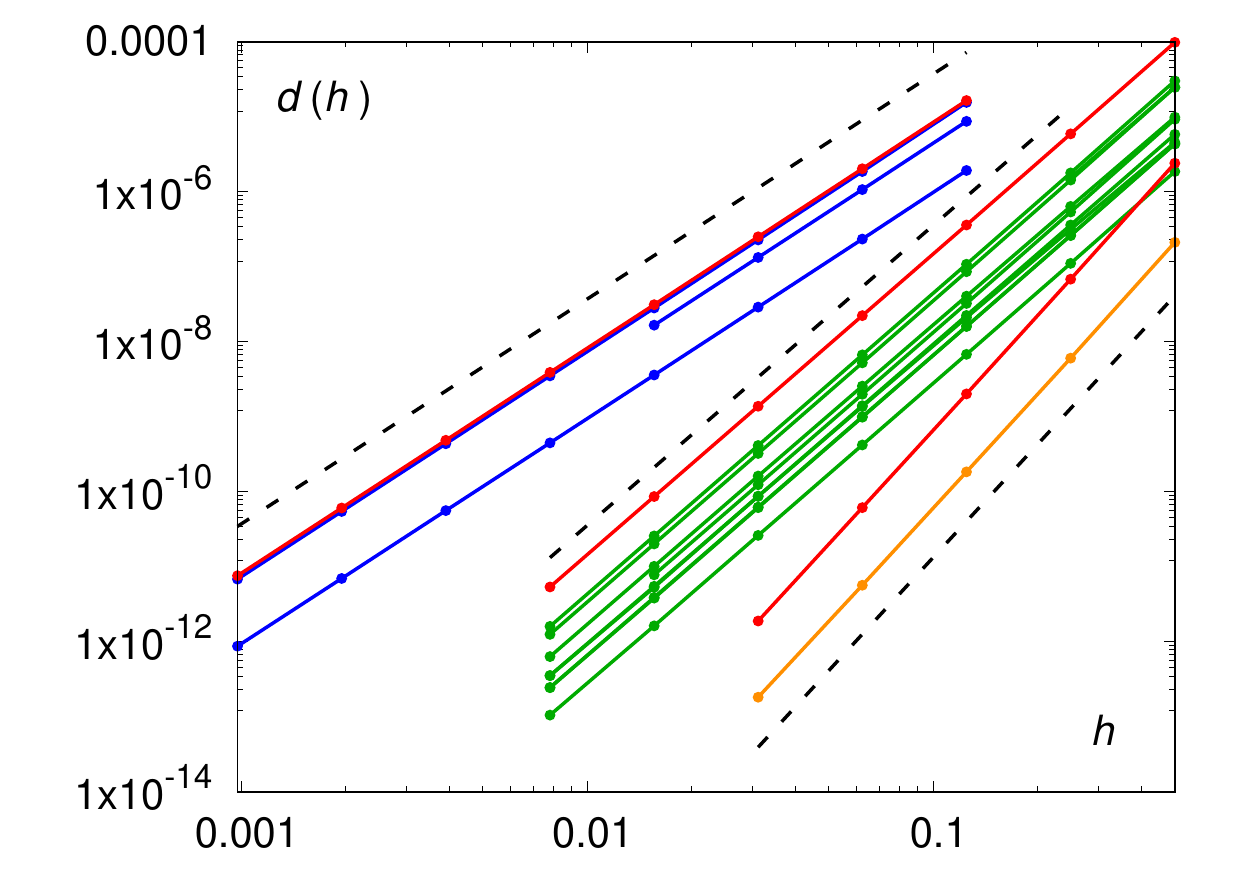}
	\caption{
		Distance from the reference solution $d(h)$ for various
		integrators as function of step size $h$ shown in a log-log plot
		for the $SO(5)$ manifold in example 2, Eq.~(\ref{eq_ex2diag}).
		As in Fig.~\ref{fig_ex1scale}, the red lines represent
		the three integrators of Munthe-Kaas type of order $p=3$, $4$ and $5$.
		The three blue lines represent the three integrators of order $p=3$ from
		Table~\ref{tab_ints}, the ten green lines the ten integrators of
		order $p=4$ and the orange line one integrator of order $p=5$ from the
		same table.
		For the $p=4$ integrators the minimum step size shown is $1/128$ and
		for $p=5$ $1/32$ since when $d(h)$ becomes comparable to $10^{-13}$
		the roundoff errors prevent correct scaling behavior. The black dashed
		lines are shown to guide the eye and represent from top to bottom
		$h^3$, $h^4$ and $h^5$, respectively.
		\label{fig_ex2scale}
	}
\end{figure}

Our second numerical example is the one\footnote{up to the dimension: we use
$5\times5$ orthogonal matrices and Ref.~\cite{MUNTHEKAAS1999115} used
$4\times4$}
used in~\cite{MUNTHEKAAS1999115},
where $Y$ is an $SO(5)$ matrix and the skew-symmetric matrix $A(Y)$
on the right hand side in Matlab notation is
\begin{equation}
\label{eq_ex2diag}
A(Y)={\rm diag}({\rm diag}(Y,+1),+1)-{\rm diag}({\rm diag}(Y,+1),-1).
\end{equation}
The inital condition $Y(0)=Y_0$ is produced randomly with
\begin{equation}
{\rm rand}(\mbox{`seed'},0);\,\,\,\,\,\,\,
[Y_0,R] = {\rm qr}({\rm rand}(5,5)).
\end{equation}

We integrate the equation of motion from $t=0$ to $t=5$ using the step
size $h=1/2^n$, $n=1,\dots,10$.
As the reference solution at time $t=5$ $Y_{ref}(t=5)$ we use 
the solution produced by the package \texttt{DiffMan}~\cite{ENGO2001323}
with the RKMK method
of order $p=6$ \texttt{butcher6} with the step size $h=1/512$.
As in the previous example we define the distance from the reference solution
\begin{equation}
d(h)=|Y(t=5,h)-Y_{ref}(t=5)|,
\end{equation}
where $|\dots|$ is the matrix 2-norm, evaluated in Matlab as \texttt{norm}.
The results for $d(h)$ are shown in Fig.~\ref{fig_ex2scale}.
The integrators again show expected scaling supporting the conjecture stated
in Sec.~\ref{sec_newconj}.

\subsection{Example 3: $SU(3)$ gradient flow}
\label{sec_ex3}

\begin{figure}[t]
	\centering
	\includegraphics[width=0.7\textwidth]{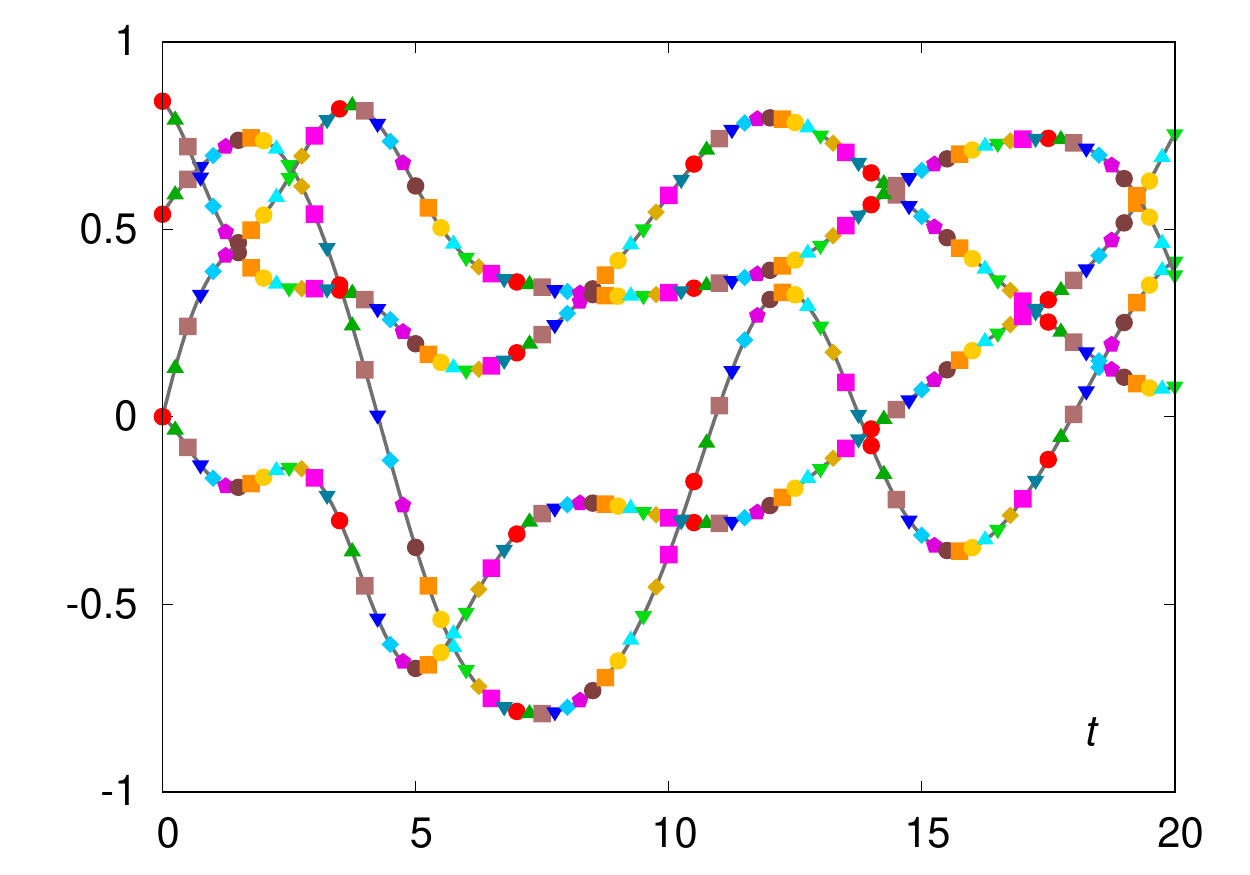}
	\caption{
		Comparison of the reference solution (see text) for the real and
		imaginary parts of the $SU(3)$ matrix elements $Y_{11}(t)$, $Y_{12}(t)$,
		shown as the four gray lines,
		with the results produced by
		the fourteen integrators listed in Table~\ref{tab_ints}.
		If plotted as lines, all results are indistinguishable from the
		exact solution. Therefore we plot the results from different
		integrators as symbols of different shape and color skipping 140 steps in
		the sequence and starting to plot 
		the first integrator at a shift of 0 steps, 
		second at a shift of 10 steps and so on.
		\label{fig_ex3sol}
	}
\end{figure}

\begin{figure}[t]
	\centering
	\includegraphics[width=0.7\textwidth]{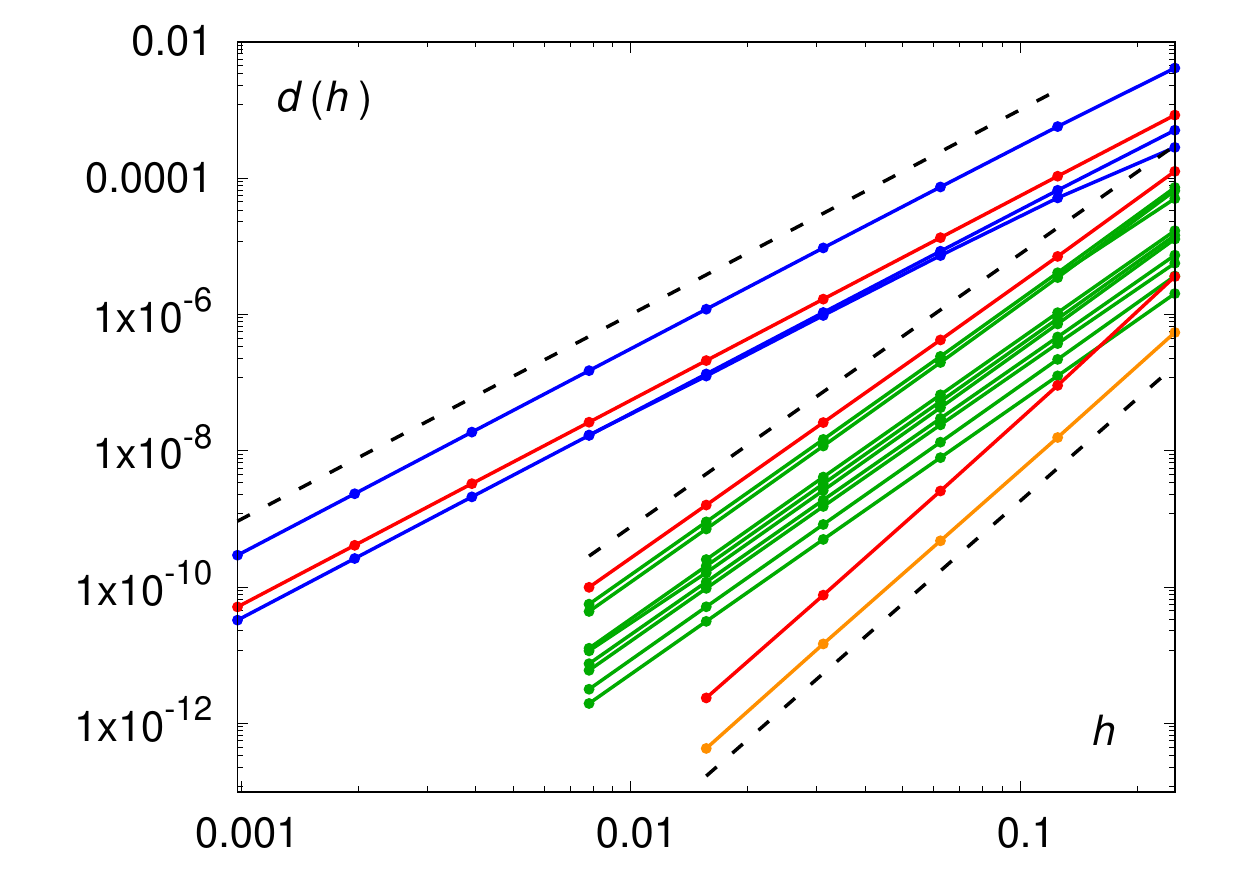}
	\caption{
		Distance from the reference solution $d(h)$ for various
		integrators as function of step size $h$ shown in a log-log plot
		for the $SU(3)$ manifold in example 3, Eq.~(\ref{eq_Ysu3}).
		As in Fig.~\ref{fig_ex1scale}, the red lines represent
		the three integrators of Munthe-Kaas type of order $p=3$, $4$ and $5$.
		The three blue lines represent the three integrators of order $p=3$ from
		Table~\ref{tab_ints}, the ten green lines the ten integrators of
		order $p=4$ and the orange line one integrator of order $p=5$ from the
		same table.
		For the $p=4$ integrators the minimum step size shown is $1/128$ and
		for $p=5$ $1/64$ since when $d(h)$ becomes comparable to $10^{-13}$
		the roundoff errors prevent correct scaling behavior. The black dashed
		lines are shown to guide the eye and represent from top to bottom
		$h^3$, $h^4$ and $h^5$, respectively.
		\label{fig_ex3scale}
	}
\end{figure}

As a third example we consider an application that is relevant for a non-perturbative
approach to quantum field theory called \textit{lattice gauge theory}.
In this case the degrees of freedom are $SU(3)$ group elements that reside on
links of a four-dimensional space-time grid and the interactions in the system
are encoded in traces of products of the $SU(3)$ matrices taken along closed
paths on the grid. L\"{u}scher in Ref.~\cite{Luscher:2010iy}
introduced a diffusion-like
procedure that suppresses short-wavelength fluctuations in the system.
This procedure leads to the following equation on the $SU(3)$ manifold:
\begin{equation}
\frac{dY}{dt}=-{\cal P}\left\{HY\right\}Y,
\label{eq_Ysu3}
\end{equation}
where $Y\in SU(3)$ and $H\in GL(3,C)$ encodes the interactions with the
neighboring degrees of freedom on the grid.
Here for numerical experiments we consider a single degree of freedom $Y$
in the presence of a fixed background $H$. The projection
\begin{equation}
{\cal P}\left\{M\right\} = \frac{1}{2}\left(M-M^\dagger\right) - 
\frac{1}{6}{\rm Tr}\left(M-M^\dagger\right)
\end{equation}
produces an element of the algebra $\mathfrak{su}(3)$ and
the right hand side of Eq.~(\ref{eq_dYAY}) in this case is
\begin{equation}
A(Y)=-{\cal P}\left\{HY\right\},
\end{equation}
where $H$ is constant. We choose $H$ as a random $3\times 3$ complex matrix
and take a diagonal initial condition $Y(t=0)=diag(e^{i},e^{i},e^{-2i})$.
Here again, to test the implementation of the integrators, we compare the
trajectory integrated with the fourteen methods of Table~\ref{tab_ints} with
the solution obtained with \texttt{DiffMan} with the RKMK integrator \texttt{butcher6}
at step size $h=1/512$. The results for the real and imaginary parts of the 
matrix elements $Y_{11}(t)$ and $Y_{12}(t)$ are shown in Fig.~\ref{fig_ex3sol}.

For the scaling study we use the same random matrix $H$ and the same initial
condition. The trajectory is integrated from $t=0$ to $t=10$
with $h=1/2^n$, $n=1,\dots,10$
and as before
we use the 2-norm as a measure of deviation from the reference solution
\begin{equation}
d(h)=|Y(t=10,h)-Y_{ref}(t=10)|.
\end{equation}
The results for $d(h)$ are shown in Fig.~\ref{fig_ex3scale}.

The integrator suggested by L\"{u}scher in the Appendix of Ref.~\cite{Luscher:2010iy}
for integrating this system, Eq.~(\ref{eq_Ysu3}) is, in fact, a 3-stage third-order
$2N$-storage integrator of the family (\ref{eq_Wc3_again}) for which we have
proven in Sec.~\ref{sec_myoc} that all integrators of this family are
third-order Lie group methods. The choice of the classical coefficients 
equivalent to the scheme in~\cite{Luscher:2010iy} is
\begin{eqnarray}
a_{21} &=& \phantom{-}\frac{1}{4},\label{eq_close_a21}\\
a_{31} &=& -\frac{2}{9},\\
a_{32} &=& \phantom{-}\frac{8}{9},\\
b_1  &=& \phantom{-}\frac{1}{4},\\
b_2  &=& \phantom{-}0,\\
b_3  &=& \phantom{-}\frac{3}{4}.\label{eq_close_b3}
\end{eqnarray}
Although the integrator in Ref.~\cite{Luscher:2010iy} was not
written in the $2N$-storage format, it was realized there that this
scheme can be used as a low-storage scheme\footnote{The simplification that
allows for this observation can be traced to the fact that $b_2=0$.},
independently of the
earlier work~\cite{WILLIAMSON198048}.
Applications to lattice gauge theory is a case where low-storage schemes
are especially attractive, since realistic systems include grids of the size
up to $96^3\times192$~\cite{Bazavov:2015yea}
which translates to about $1.2\times10^{10}$ double precision numbers to be
stored on a (super)computer just to represent the system. 
While using a $2N$-storage scheme requires twice that amount, an
equivalent 3-stage third-order RKMK method requires four times this amount,
with further increasing requirements for higher order schemes.
Since the conjecture about higher than order $p=3$ $2N$-storage classical RK
methods holds true numerically for the $SU(3)$ case,
as illustrated in Fig.~\ref{fig_ex3scale}, this opens a possibility of constructing
higher order Lie group $2N$-storage methods for applications in
lattice gauge theory.

\subsection{Example 4: Van der Pol oscillator}

\begin{figure}[t]
	\centering
	\includegraphics[width=0.7\textwidth]{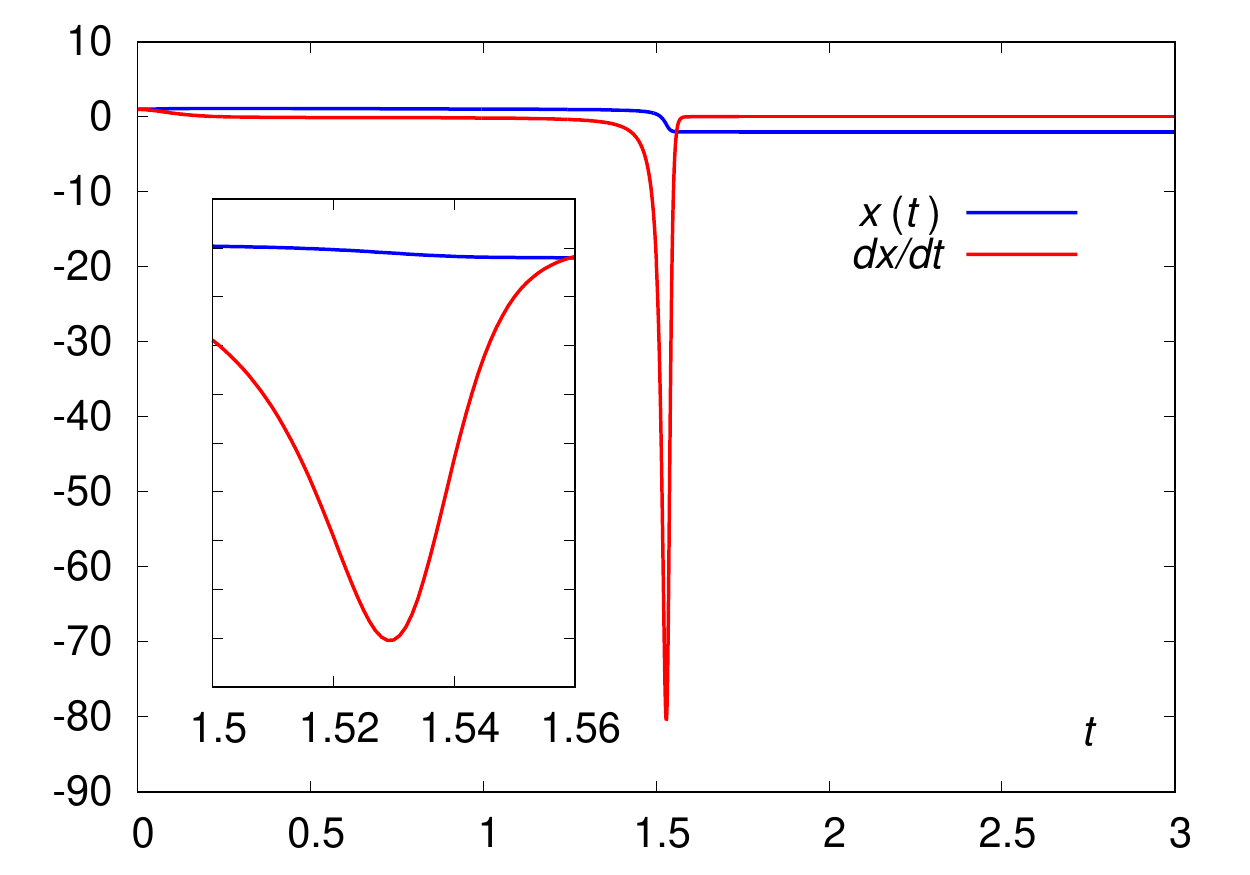}
	\caption{
		The dependence of the coordinate $x(t)$ (blue) and
		the velocity $\dot{x}(t)$ (red) on time for the van der Pol system,
		Eq.~(\ref{eq_vdp}). This solution is
		produced with the BWRRK33 integrator with step size $h=0.001$. The other
		integrators produce results indistinguishable on the scale of the figure
		and are not shown. The inset magnifies the 
		horizontal scale in the vicinity of 
		the ``needle''.
		\label{fig_ex4sol}
	}
\end{figure}

\begin{figure}[t]
	\centering
	\includegraphics[width=0.7\textwidth]{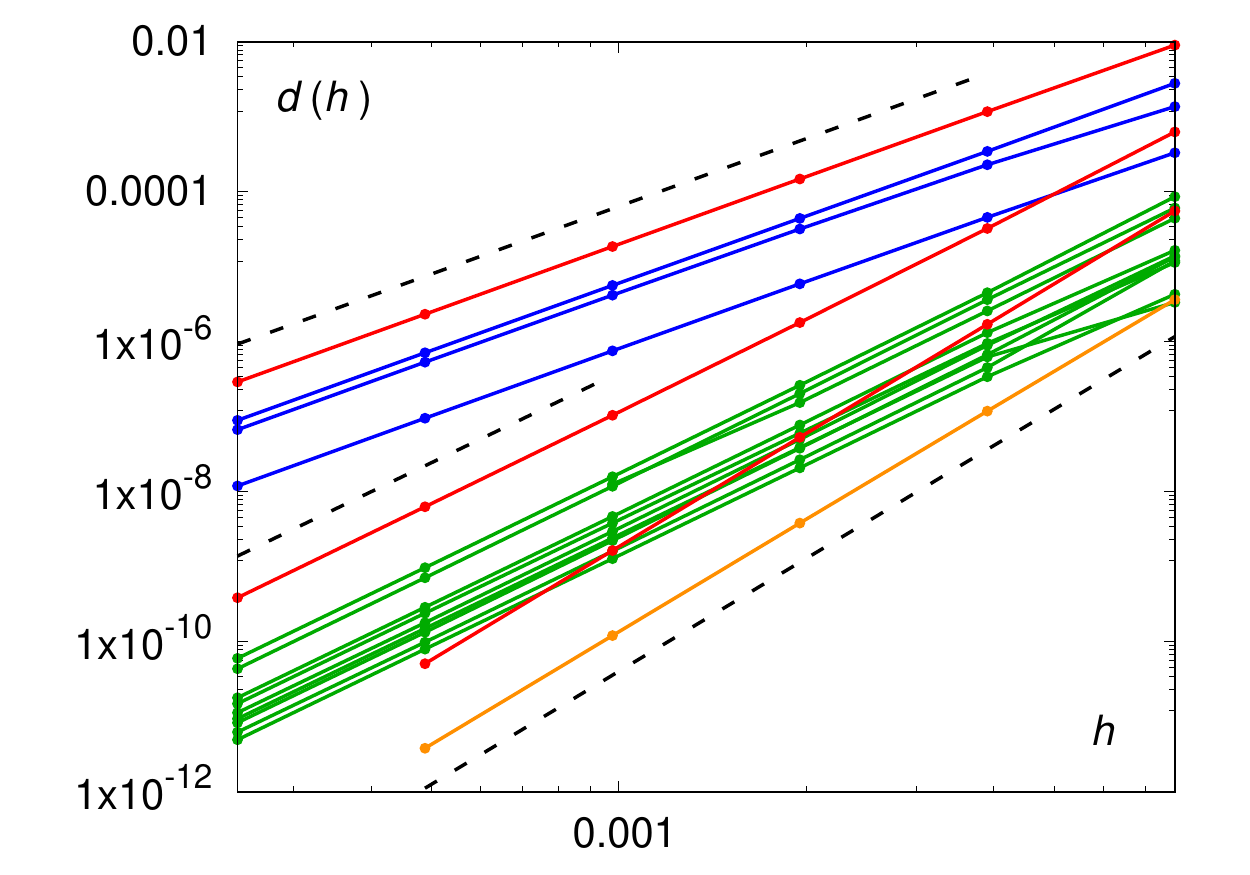}
	\caption{
		Distance from the reference solution $d(h)$ for various
		integrators as function of step size $h$ shown in a log-log plot
		for the van der Pol oscillator in example 4, Eq.~(\ref{eq_vdp}).
		As in Fig.~\ref{fig_ex1scale}, the red lines represent
		the three integrators of Munthe-Kaas type of order $p=3$, $4$ and $5$.
		The three blue lines represent the three integrators of order $p=3$ from
		Table~\ref{tab_ints}, the ten green lines the ten integrators of
		order $p=4$ and the orange line one integrator of order $p=5$ from the
		same table.
		For the $p=5$ integrators the minimum step size shown is $1/2048$
		since when $d(h)$ becomes comparable to $10^{-13}$
		the roundoff errors prevent correct scaling behavior. The black dashed
		lines are shown to guide the eye and represent from top to bottom
		$h^3$, $h^4$ and $h^5$, respectively.
		\label{fig_ex4scale}
	}
\end{figure}

The van der Pol equation
\begin{equation}
\label{eq_vdp}
\frac{d^2x}{dt^2}-\mu(1-x^2)\frac{dx}{dt}+x=0
\end{equation}
has also been used as a test case in the literature
on Lie group methods~\cite{Curry2019}.
In this case, a Lie group method is used as an
exponential integrator that may handle stiff systems better than classical
RK schemes. Eq.~(\ref{eq_vdp}) can be
rewritten in a vector form
\begin{equation}
\frac{d}{dt}\left(
\begin{array}{c}
x \\
\dot{x}
\end{array}
\right)=
\left(
\begin{array}{rc}
0 & 1 \\
-1 & \mu(1-x^2)
\end{array}
\right)
\left(
\begin{array}{c}
x \\
\dot{x}
\end{array}
\right)
\end{equation}
where we can identify $Y$ as a two-dimensional vector and $A(Y)\in GL(2,R)$.
As in Ref.~\cite{Curry2019} we choose $\mu=60$ and the initial condition
\begin{equation}
Y(0)=\left(
\begin{array}{c}
1 \\
1
\end{array}
\right).
\end{equation}
At such a large value of $\mu$ the system is stiff as shown in Fig.~\ref{fig_ex4sol}
and the ``needle'' occurs approximately at $t=1.53$.

We integrate the system from $t=0$ to $t=2$, \textit{i.e.} past the ``needle'',
with step size $h=1/2^n$, $n=7,\dots,12$.
As a reference solution $Y_{ref}(t=2)$ we use the result from the \texttt{DiffMan}
package with the RKMK integrator \texttt{butcher6} at step size $h=1/4096$ and
define the distance from the reference solution $d(h)$ the same way as in the
Example 1 in Sec.~\ref{sec_ex1}. The results for $d(h)$ for various 
$2N$-storage schemes of Table~\ref{tab_ints} are shown in Fig.~\ref{fig_ex4scale}.

Notice that given the stiffness of the system, we use, in general, a range of
smaller step sizes than in the other examples, but all the integrators do show
the scaling expected from the conjecture in Sec.~\ref{sec_newconj}.

\subsection{Example 5: Non-autonomous problem in $SO(3)$}

\begin{figure}[t]
	\centering
	\includegraphics[width=0.7\textwidth]{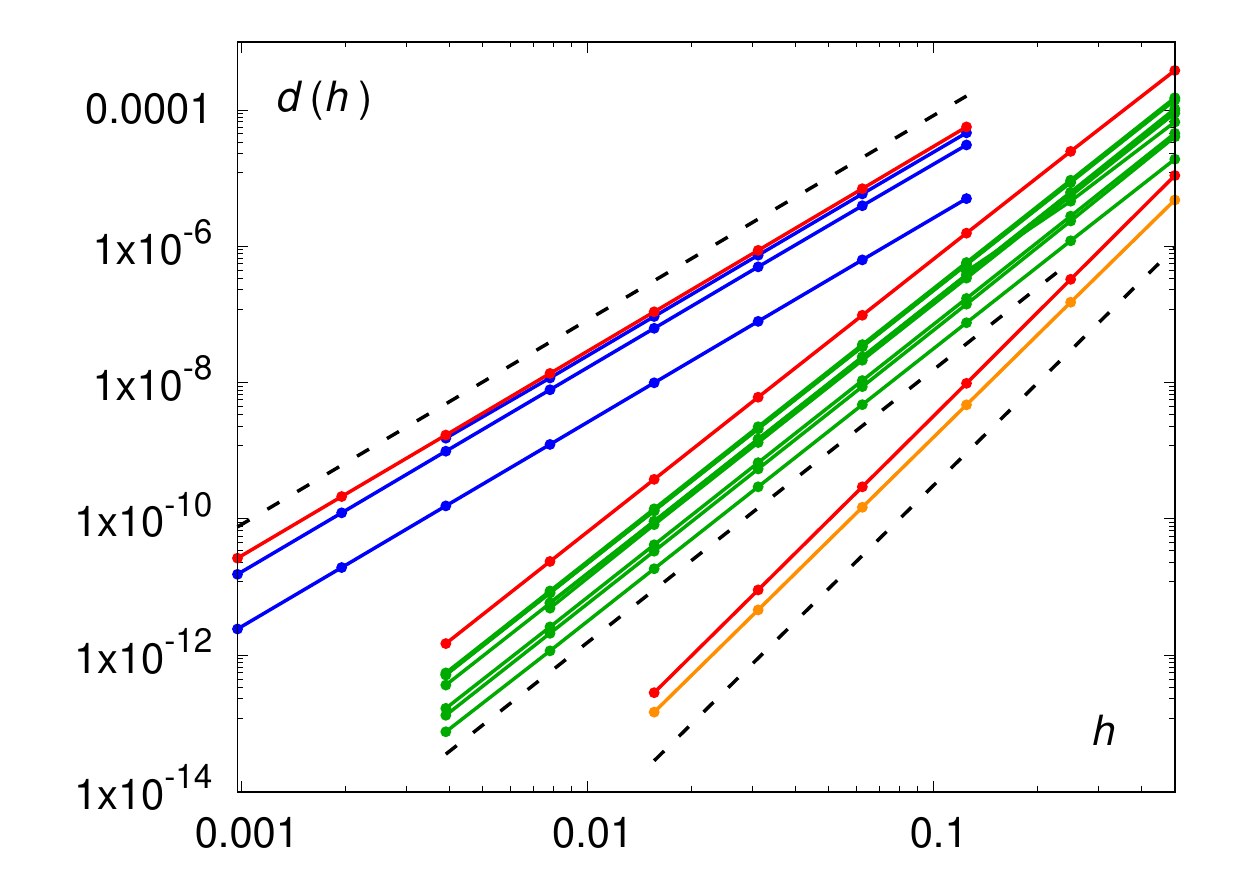}
	\caption{
		Distance from the reference solution $d(h)$ for various
		integrators as function of step size $h$ shown in a log-log plot
		for the $SO(3)$ matrix $Y$ in example 5
		(non-autonomous problem), Eq.~(\ref{eq_AtY}).
		As in Fig.~\ref{fig_ex1scale}, the red lines represent
		the three integrators of Munthe-Kaas type of order $p=3$, $4$ and $5$.
		The three blue lines represent the three integrators of order $p=3$ from
		Table~\ref{tab_ints}, the ten green lines the ten integrators of
		order $p=4$ and the orange line one integrator of order $p=5$ from the
		same table.
		For the $p=4$ integrators the minimum step size shown is $1/256$ and
		for $p=5$ $1/64$ since when $d(h)$ becomes comparable to $10^{-13}$
		the roundoff errors prevent correct scaling behavior. The black dashed
		lines are shown to guide the eye and represent from top to bottom
		$h^3$, $h^4$ and $h^5$, respectively.
		\label{fig_ex5scale}
	}
\end{figure}

As the final example we consider a non-autonomous problem which is included as
one of the examples in \texttt{DiffMan}~\cite{ENGO2001323}: $Y\in SO(3)$ and
\begin{equation}
\label{eq_AtY}
A(t,Y)=
\left(
\begin{array}{rrr}
0 & t\phantom{^2} & 1\phantom{^2} \\
-t & 0\phantom{^2} & -t^2 \\
-1 & t^2 & 0\phantom{^2}
\end{array}
\right)\in\mathfrak{so}(3).
\end{equation}
This test is different from the previous ones since the coefficients $c_i$,
Eq.~(\ref{eq_ci}), now 
enter the game and we investigate if that may lead to a breakdown of the
conjecture of Sec.~\ref{sec_newconj}. We choose a unit matrix as the initial condition,
as in \texttt{DiffMan}, and integrate the trajectory from $t=0$ to $t=1$ with
step sizes $h=1/2^n$, $n=1,\dots,10$. As the reference solution we take the result
from \texttt{DiffMan} integrated with \texttt{butcher6} with $h=1/1024$.
The distance from the reference solution $d(h)$ is defined the same way as in
examples 2, Sec.~\ref{sec_ex2} and 3, Sec.~\ref{sec_ex3}. The results for $d(h)$
for the methods of Table~\ref{tab_ints} is shown in Fig.~\ref{fig_ex5scale}.
As can be observed from the figure, the integrators again show the expected
scaling.

\section{Conclusions}
\label{sec_concl}
We have shown in Sec.~\ref{sec_myoc} that 3-stage third-order classical $2N$-storage
Runge-Kutta methods of Ref.~\cite{WILLIAMSON198048} are also third-order 
commutator-free Lie group methods, since the coefficients satisfy the same order conditions
in both cases: four classical ones and an additional one arising either from writing
the scheme in $2N$-storage format~\cite{WILLIAMSON198048} or constructing a commutator-free
Lie group integrator in a specific format proposed in 
Eqs.~(\ref{eq_exp3_y1})--(\ref{eq_exp3_yth}) that automatically leads to
Eqs.~(\ref{eq_myY1})--(\ref{eq_myYth}).

Given this similarity, we conjectured in Sec.~\ref{sec_newconj} that this observation
holds also for third-order $2N$-storage 
classical RK methods with more than three stages
and for methods of order four and five. In Sec.~\ref{sec_num} we considered
five different numerical examples and studied how the global error of the $2N$-storage
classical RK methods available in the 
literature~\cite{CK1994,BERNARDINI20094182,BERLAND20061459,TOULORGE20122067,STANESCU1998674,ALLAMPALLI20093837,NIEGEMANN2012364,Yan2017}
scales with the step size when
the schemes are used as commutator-free Lie group methods of the $2N$-storage
format proposed in Eqs.~(\ref{eq_mydy})--(\ref{eq_myY}). We found that in all test
cases numerical evidence supports the conjecture.

As a next step, it is obviously desirable
to find an analytic proof of the conjecture
of Sec.~\ref{sec_newconj}.
In the meantime, one can check the order conditions for a Lie group integrator
based on a given $2N$-storage coefficients scheme with the methods of
Ref.~\cite{Knoth2020git} or determine what order is achieved numerically
as in the examples presented in Sec.~\ref{sec_num}.

If the conjecture is proven to be correct, there are two possible benefits of the
proposed low-storage commutator-free Lie group integrators. First, in large-scale
calculations one can significantly reduce memory requirements compared to 
the available Lie group methods (both, commutator-free and with commutators) and
also reuse the exponentials at every step of the calculation which leads to the
requirement of evaluating exactly $s$ exponentials for a $s$-stage method.
Second, it may be easier to develop new schemes of this type for differential
equations on manifolds in a way it has been done
for the classical $2N$-storage methods. Once the scheme is written in a $2N$-storage
format, Eqs.~(\ref{eq_mydy})--(\ref{eq_myY}),
one needs to find the coefficients of the $2N$-storage scheme from
satisfying \textit{only} \textit{classical} Runge-Kutta order conditions. The 
property that the scheme is
a $2N$-storage scheme will automatically (again, if the conjecture is true) 
satisfy all the additional constraints arising from non-commutativity.

\bigskip
\textbf{Acknowledgements.}
I thank Andrea Shindler for careful reading
and comments on the manuscript and Oswald Knoth for bringing my attention
to Ref.~\cite{KNOTH1998327,Wensch2009}, comments on the manuscript and,
most importantly, for independently
checking the order conditions for all coefficient schemes of
Table~\ref{tab_ints} (as Lie group integrators of the form 
(\ref{eq_mydy})--(\ref{eq_myY}))
with his software available at~\cite{Knoth2020git}
and explaining me the theory and algorithmic details behind that
software.
This work was in part supported
by the U.S. National Science Foundation under award
PHY-1812332.

\appendix

\section{Matlab code}
\label{app_matlab}
To illustrate the usage of the $2N$-storage schemes as commutator-free
Lie group integrators we present below a Matlab script that generates
a scaling plot similar to Fig.~\ref{fig_ex1scale} using the integrators
BWRRK33, TSRKF84 and YRK135 from Table~\ref{tab_ints}.

\begin{verbatim}
%%%%%%%%%%%%%%%%%%%%%%%%%%%%%%%%
% main function in the script:
% get scaling with step size
% by comparing to exact solution
%%%%%%%%%%%%%%%%%%%%%%%%%%%%%%%%
function lie_integrate

format long;

global I;

% initial condition
% second component is 0 to match
% simplifying assumptions of exact solution
% in Marsden, Ratiu, Introduction to Mechanics and Symmetry
y0 = [ -sqrt(8)/3; 0; 1/3 ];

% moment of intertia, same as in
% Celledoni, Marthinsen, Owren
% Future Generation Computer Systems, 19 (2003) 341-352
I1 = 7/8;
I2 = 5/8;
I3 = 1/4;
I = diag( [ I1 I2 I3 ] );

% integration parameters
T0 = 0;
T = 3;

% set array with time steps
array_dt = ...
[ 1/2048 1/1024 1/512 1/256 1/128 1/64 1/32 1/16 1/8 1/4 1/2 ];
Ndt = length( array_dt );

% number of integrators
Nint = 3;

% storage
error_dt = zeros( Ndt, Nint );
sol = zeros( 3, Ndt, Nint);

% coefficients for low-storage integrators
% 3-stage third-order BWRRK33
Nstages_3 = 3;
A_3 = [ 0 ...
       -0.637694471842202 ...
       -1.306647717737108 ];
B_3 = [ 0.457379997569388 ...
        0.925296410920922 ...
        0.393813594675071 ];
C_3 = [ 0 ...
        0.457379997569388 ...
        0.792620002430607 ];
% 8-stage fourth-order TSRKF84 of
% Toulorge, Desmet, J. Comp. Phys. 231 (2012) 2067-2091
Nstages_4 = 8;
A_4 = [ 0 ...
       -0.5534431294501569 ...
        0.01065987570203490 ...
       -0.5515812888932000 ...
       -1.885790377558741 ...
       -5.701295742793264 ...
        2.113903965664793 ...
       -0.5339578826675280 ];
B_4 = [ 0.08037936882736950 ...
        0.5388497458569843 ...
        0.01974974409031960 ...
        0.09911841297339970 ...
        0.7466920411064123 ...
        1.679584245618894 ...
        0.2433728067008188 ...
        0.1422730459001373 ];
C_4 = [ 0 ...
        0.08037936882736950 ...
        0.3210064250338430 ...
        0.3408501826604660 ...
        0.3850364824285470 ...
        0.5040052477534100 ...
        0.6578977561168540 ...
        0.9484087623348481 ];
% 13-stage fifth-order YRK135 of
% Yan, Chin. J. Chem. Phys 30 (2017) 277-286
Nstages_5 = 13;
A_5 = [ 0 ...
       -0.33672143119427413 ...
       -1.2018205782908164 ...
       -2.6261919625495068 ...
       -1.5418507843260567 ...
       -0.2845614242371758 ...
       -0.1700096844304301 ...
       -1.0839412680446804 ...
       -11.61787957751822 ...
       -4.5205208057464192 ...
       -35.86177355832474 ...
       -0.000021340899996007288 ...
       -0.066311516687861348 ];
B_5 = [ 0.069632640247059393 ...
        0.088918462778092020 ...
        1.0461490123426779 ...
        0.42761794305080487 ...
        0.20975844551667144 ...
       -0.11457151862012136 ...
       -0.01392019988507068 ...
        4.0330655626956709 ...
        0.35106846752457162 ...
       -0.16066651367556576 ...
       -0.0058633163225038929 ...
        0.077296133865151863 ...
        0.054301254676908338 ];
C_5 = [ 0 ...
        0.069632640247059393 ...
        0.12861035097891748 ...
        0.34083022189561149 ...
        0.54063706308495402 ...
        0.59927749518613931 ...
        0.49382042519248519 ...
        0.48207852767699775 ...
        0.82762865209834452 ...
        0.82923953914857933 ...
        0.67190565554748019 ...
        0.87194975193167848 ...
        0.94930216564503562 ];


% parameters for exact solution
I1y0 = (I^-1) * y0;
h = 1/2 * ( y0.'*I1y0 );
y02 = y0.'*y0;
a = y02 / (2*h);
b = 2*h / sqrt(y02);
alpha = sqrt( a*I2*(a-I3)/(I2-I3) ) * b;
beta = sqrt( a*I2*(I1-a)/(I1-I2) ) * b;
mu = sqrt( a*(I1-a)*(I2-I3)/(I1*I2*I3) ) * b;
k = sqrt( (I1-I2)*(a-I3)/(I1-a)/(I2-I3) );
delta = sqrt( I3*(I1-a)*a/(I1-I3) ) * b;
gamma = sqrt( I1*(a-I3)*a/(I1-I3) ) * b;

% exact solution at T
snmt = jacobiSN( mu*T, k^2 );
cnmt = jacobiCN( mu*T, k^2 );
dnmt = jacobiDN( mu*T, k^2 );
y_exact = [ -gamma * cnmt; alpha * snmt; delta * dnmt ];


% loop over step sizes
for idt = 1:Ndt

 % integrate with current step size
 for iint=1:Nint
  dt = array_dt(idt);
  if iint==1
   sol( :, idt, iint ) = ...
   integrate( T0, T, dt, y0, @rhs_f, Nstages_3, A_3, B_3, C_3 );
  elseif iint==2
   sol( :, idt, iint ) = ...
   integrate( T0, T, dt, y0, @rhs_f, Nstages_4, A_4, B_4, C_4 );
  elseif iint==3
   sol( :, idt, iint ) = ...
   integrate( T0, T, dt, y0, @rhs_f, Nstages_5, A_5, B_5, C_5 );
  end

  % get error
  error_dt( idt, iint ) = norm( sol( :, idt, iint ) - y_exact );
 end

end

% plot
for iint=1:Nint
 if iint==1
  cc = [1 0 0];
  xaxis = array_dt;
  yaxis = error_dt( :, iint );
 elseif iint==2
  cc = [0 0.7 0];
  xaxis = array_dt( 3:Ndt );
  yaxis = error_dt( 3:Ndt, iint );
 elseif iint==3
  cc = [0 0 1];
  xaxis = array_dt( 5:Ndt );
  yaxis = error_dt( 5:Ndt, iint );
 end
 loglog( xaxis, yaxis, 'Color', cc, 'LineWidth', 3 );
 if iint==1
  hold on
 end
end
hold off

end

%%%%%%%%%%%%%%%%%%%%%%%%%%%%
% right hand side function
% t - current time
% y - current function value
%%%%%%%%%%%%%%%%%%%%%%%%%%%%
function r = rhs_f( t, y )

 global I;

 I1y = (I^-1) * y;
 r = zeros(3,3);
 r(1,2) =  I1y(3);
 r(2,1) = -I1y(3);
 r(1,3) = -I1y(2);
 r(3,1) =  I1y(2);
 r(2,3) =  I1y(1);
 r(3,2) = -I1y(1);

end

%%%%%%%%%%%%%%%%%%%%%%%%%%%%%%%%%%%%%%%%%%%%%%%%%%
% low-storage commutator-free Lie group integrator
% T0 - initial time
% T - final time
% dt - step size
% Y0 - initial function value
% rhs - function to evaluate right hand side
% Ns - number of stages of the integrator
% A, B, C - coefficients in 2N-storage format
% Note: A(1)=0 is required
% result - function value at time T Y(t=T)
%%%%%%%%%%%%%%%%%%%%%%%%%%%%%%%%%%%%%%%%%%%%%%%%%%
function result = integrate( T0, T, dt, Y0, rhs, Ns, A, B, C )

 % current time
 Tcur = T0;
 % initial function value
 Y = Y0;
 dY = 0;

 while Tcur < T
  if dt > T-Tcur
   dt = T - Tcur;
  end

  for k=1:Ns
   dY = A(k)*dY + dt*rhs( Tcur + C(k)*dt, Y );
   Y = expm( B(k)*dY )*Y;
  end

  % set values for next iteration
  Tcur = Tcur + dt;
 end

 result = Y;

end

\end{verbatim}





\bibliography{lie_int}




\end{document}